\documentclass[11pt]{article}
\usepackage{amssymb,amsmath,amsthm,subfigure,url}
\usepackage{mathrsfs}
\usepackage[pdftex]{graphicx}

\newtheorem{thm}{Theorem}[section]
\newtheorem{lemma}[thm]{Lemma}
\newtheorem{prop}[thm]{Proposition}
\newtheorem{cor}[thm]{Corollary}

\theoremstyle{definition}
\newtheorem{defn}[thm]{Definition}
\newtheorem{example}[thm]{Example}
\newtheorem{construction}[thm]{Construction}

\newtheorem*{notation}{Notation}

\newcommand{\bl}{\boldsymbol{\ell}}

\renewcommand{\aa}{\mathbf{a}}
\newcommand{\bb}{\mathbf{b}}
\newcommand{\kk}{\mathbf{k}}
\newcommand{\rr}{\mathbf{r}}
\renewcommand{\tt}{\mathbf{t}}
\newcommand{\uu}{\mathbf{u}}
\newcommand{\vv}{\mathbf{v}}
\newcommand{\ww}{\mathbf{w}}
\newcommand{\xx}{\mathbf{x}}
\newcommand{\yy}{\mathbf{y}}

\renewcommand{\AA}{\mathbf{A}}
\newcommand{\BB}{\mathbf{B}}
\newcommand{\CC}{\mathbf{C}}
\newcommand{\DD}{\mathbf{D}}
\newcommand{\RR}{\mathbf{R}}
\renewcommand{\SS}{\mathbf{S}}
\newcommand{\TT}{\mathbf{T}}
\newcommand{\XX}{\mathbf{X}}

\newcommand{\vmax}{v_\textrm{max}}
\newcommand{\vvmax}{\mathbf{v}_\textrm{max}}
\newcommand{\kmin}{k_\textrm{min}}
\newcommand{\kkmin}{\mathbf{k}_\textrm{min}}

\renewcommand{\hat}{\widehat}

\newcommand\ceilfrac[2]{ \left\lceil \frac{#1}{#2} \right\rceil }
\DeclareMathOperator{\cat}{cat}

\renewcommand{\emptyset}{\varnothing}

\title{Generalized covering designs and clique coverings}% \\REVISED VERSION--1st DRAFT}
\author{%
Robert F.~Bailey\footnote{Corresponding author.},\footnote{Department of Mathematics and Statistics, University of Regina, 3737 Wascana Parkway, Regina, SK, S4S 0A2, Canada.  E-mail: \texttt{robert.bailey@uregina.ca}}\,\,
Andrea C.~Burgess\footnote{Department of Mathematics, Ryerson University, 350 Victoria St., Toronto, ON, M5B 2K3, Canada.  E-mail: \texttt{andrea.burgess@ryerson.ca}},\,\,
Michael S.~Cavers\footnote{Department of Mathematics and Statistics, University of Calgary, 2500 University Drive NW, Calgary, AB, T2N 1N4, Canada.  E-mail: \texttt{michael.cavers@ucalgary.ca}},\,\,
Karen Meagher\footnote{Department of Mathematics and Statistics, University of Regina, 3737 Wascana Parkway, Regina, SK, S4S 0A2, Canada.  E-mail: \texttt{karen.meagher@uregina.ca}}
}

\begin{document}

\maketitle

\begin{abstract}
Inspired by the ``generalized $t$-designs'' defined by Cameron [P.~J.~Cameron, A generalisation of $t$-designs, {\em Discrete Math.} {\bf 309} (2009), 4835--4842], we define a new class of combinatorial designs which simultaneously provide a generalization of both covering designs and covering arrays.  We then obtain a number of bounds on the minimum sizes of these designs, and describe some methods of constructing them, which in some cases we prove are optimal.  Many of our results are obtained from an interpretation of these designs in terms of clique coverings of graphs.\\


\end{abstract}

\section{Introduction} \label{section:intro}
In a 2009 paper \cite{Cameron09}, Cameron introduced a new class of combinatorial designs, which simultaneously generalizes various well-known classes of designs, including $t$-designs, mutually orthogonal Latin squares, orthogonal arrays and 1-factorizations of complete graphs.  %%\marginpar{\sf [ADDITION]} \\
Further work on Cameron's ``generalized $t$-designs'' has been done by Soicher \cite{Soicher} and others \cite{Drizen,Patterson}, while the earlier papers of Martin \cite{Martin98,Martin99} and Teirlinck \cite{Teirlinck89} discuss related objects.  
In a remark near the end of his paper, Cameron suggests that a similar definition can be made for generalizing covering designs.

The purpose of this paper is to pursue such a generalization, i.e.\ to define a broad class of combinatorial designs with a ``covering'' property which includes previously well-studied, and widely-applied, families of designs as special cases.  Analogous to $t$-designs and orthogonal arrays are covering designs and covering arrays, respectively; the designs we define in this paper form a simultaneous generalization of both of these ``covering'' objects. %%\marginpar{\sf [ADDITION]} %\\

The key difference when studying covering problems rather than ``ordinary'' designs is that the question is typically not whether the designs exist (this is usually trivial to answer), but obtaining bounds on the minimum size, and constructing optimal (or near-optimal) designs.

Background material on most classes of designs can be found in the {\em Handbook of Combinatorial Designs} \cite{handbook}.  However, we give some relevant definitions here.

\subsection{Covering designs} \label{subsect:covdes}

\begin{defn} \label{defn:covdes}
Let $v,k,t,\lambda$ be positive integers with $v\geq k\geq t$.  A {\em $(v,k,t)_\lambda$-covering design} is a family $\mathcal{C}$ of $k$-subsets (called {\em blocks}) of a $v$-set $X$, where any $t$-subset of $X$ is contained in at least $\lambda$ members of $\mathcal{C}$.
\end{defn}

%%\marginpar{\sf [ADDITION]}
The notation ``$t$-$(v,k,\lambda)$ covering design'' is also used in the literature (in \cite{handbook,MillsMullin92}, for example).
In the case where each $t$-subset occurs {\em exactly} $\lambda$ times, we have a {\em $t$-$(v,k,\lambda)$ design}.  We remark that the definition allows for blocks to be repeated.  A survey of results on covering designs can be found in Mills and Mullin \cite{MillsMullin92}.  Usually, we are only concerned with the case $\lambda=1$, and omit the subscript $\lambda$.

\begin{example}
The following is an example of an $(8,5,2)$-covering design, where $X=\{1,\ldots,8\}$:
\[ \begin{array}{c}
1 \,\,\, 2 \,\,\, 3 \,\,\, 4 \,\,\, 5 \\
1 \,\,\, 5 \,\,\, 6 \,\,\, 7 \,\,\, 8 \\
2 \,\,\, 3 \,\,\, 6 \,\,\, 7 \,\,\, 8 \\
4 \,\,\, 5 \,\,\, 6 \,\,\, 7 \,\,\, 8 
\end{array} \]
By hand, it is easy to verify that each pair chosen from $X$ is contained in at least one of the 5-sets given.
\end{example}

%{\sf [BOUNDS ON $C(v,k,t)$; SCHONHEIM, ETC.]}
Unlike $t$-designs, it is clear that if $v\geq k\geq t$, then a $(v,k,t)$-covering design always exists: we simply take all the $k$-subsets of $X$ and discard any that are unnecessary.  Of course, this is not an efficient approach.  We would like our designs to have the smallest number of blocks possible: this number is called the {\em covering number} and is denoted $C(v,k,t)$.  If we have parameters $v,k,t$ for which no covering design can exist, then we say $C(v,k,t)=0$; clearly, this can only happen if $v<k$ or $k<t$.

There are many results on finding or bounding covering numbers, which can be found in \cite[section VI.11]{handbook}.  The most general bound is known as the {\em Sch\"onheim bound} \cite{Schonheim}, and is given below.

\begin{thm}[The Sch\"onheim bound] \label{thm:Schonheim}
Where $v\geq k\geq t$, we have
\[ C(v,k,t) \geq \left\lceil \frac{v}{k} \left\lceil \frac{v-1}{k-1} \cdots \left\lceil \frac{v-t+1}{k-t+1} \right\rceil \cdots \right\rceil \right\rceil. \]
\end{thm}

However, many more specific results are known, and exact values of $C(v,k,t)$ have been determined in many cases.  (Often, these are the result of sophisticated computer searches: see \cite{lajolla} for an online database of the best-known covering designs for small values of $v$, $k$ and $t$.)  In particular, in the case where $t=2$ and the ratio $v/k \leq 13/4$, the exact values of $C(v,k,t)$, and constructions of covering designs of those sizes, are all known.  For $v/k \leq 3$, the constructions are due to various authors, and can be found in a paper by Mills \cite{Mills79}.  For $3 < v/k \leq 13/4$, the constructions are much newer, and are due to Greig, Li and van Rees \cite{Greig06}.  (See \cite[Theorem VI.11.31]{handbook}, for a summary.)

%%\marginpar{\sf [ADDITION]}
An important asymptotic result is due to R\"odl \cite{Rodl85}, who in 1985 proved a conjecture of Erd{\H{o}}s and Hanani \cite{Erdos63} which asserts that for fixed values of $k$ and $t$,
\[ \lim_{v\to\infty} C(v,k,t) \frac{ {k \choose t} }{ {v \choose t} } = 1. \]
An alternative proof of R\"odl's theorem was subsequently obtained by Spencer \cite{Spencer85}, while constructions for covering designs which (asymptotically) meet this bound were obtained by Gordon {\em et al.}\ \cite{GordonKuperbergPatashnik95,Gordon96}.

Covering designs, or objects obtained from them, have been used in various applications, many of which are related to computing or communications.  These include quorum systems in distributed databases \cite{quorum}, threshold schemes in cryptography \cite{Rees99}, and decoding algorithms for error-correcting codes \cite{ecpg,Gordon82,Kroll08}. 

\subsection{Covering arrays} \label{subsect:covarray}

\begin{defn} \label{defn:covarray}
Let $N,k,s,t,\lambda$ be positive integers.  A {\em covering array} ${\rm CA}_\lambda(N;k,s,t)$ is an $N\times k$ array with entries from an alphabet of size $s$, with the property that in every set of $t$ columns, each $t$-tuple of symbols from the alphabet occurs in at least $\lambda$ rows.
\end{defn}

Note that such an array where every $t$-tuple occurs in {\em exactly} $\lambda$ rows is known as an {\em orthogonal array} (see \cite{Hedayat99}).  As with covering designs, we usually only treat the case where $\lambda=1$ and omit the subscript $\lambda$.  The parameter $t$ is called the {\em strength} of the covering array.  Usually we fix the parameters $k$, $s$ and $t$ and want to find the smallest $N$ such that there exists a ${\rm CA}(N;k,s,t)$: this value of $N$ is called the {\em covering array number}, and is denoted by ${\rm CAN}(k,s,t)$.

\begin{example}
The following is an example of a ${\rm CA}(5;4,2,2)$ where the alphabet is $\{0,1\}$:
\[ \begin{array}{cccc}
0 & 0 & 0 & 0 \\
1 & 1 & 1 & 0 \\
1 & 1 & 0 & 1 \\
1 & 0 & 1 & 1 \\
0 & 1 & 1 & 1
\end{array} \]
In each pair of columns, each of the $2^2$ possible combinations $00,01,10,11$ appears at least once.
\end{example}

Note the slight abuse of notation: in the case of covering designs, each row of the table is a set, whereas in the case of covering arrays, each row of the table is a $k$-tuple.  This distinction will be important later in the paper.

A more general object is the {\em mixed covering array}, as defined by Moura {\em et al.}\ \cite{Moura03}.  This is a covering array where each column has its own alphabet, and these may have different sizes.  If the alphabet sizes are $\vv=(v_1,\ldots,v_k)$, then such an object is denoted by ${\rm MCA}(N;k,\vv,t)$, while the least possible size is denoted by ${\rm MCAN}(k,\vv,t)$.  

Covering arrays, and mixed covering arrays, have a number of applications, most notably in software testing, but also in other areas such as  computational biology.  The survey by Colbourn \cite{Colbourn04} describes many of these applications.

\section{Generalized covering designs} \label{section:defn}

\subsection{Notation}
%{\sf [NEW SECTION ON NOTATION]}
Throughout the remainder of this paper, we will be using $m$-tuples, or vectors, of both integers and finite sets.  We begin with explaining our terminology and notation for this.

Suppose $\xx=(x_1,x_2,\ldots,x_m)$ and $\yy=(y_1,y_2,\ldots,y_m)$ are $m$-tuples of integers.  We write $\xx\leq\yy$ to mean that $x_i\leq y_i$ for all $i\in\{1,2,\ldots,m\}$.  We use the notation $x^m$ to denote the vector with $m$ entries, all equal to $x$.

There are various operations we can perform on vectors of integers, as well as the usual operations of addition, subtraction and scalar multiplication.
%Restriction
First, given a vector $\xx$ and subset $I$ of the indices $\{1,2,\ldots,m\}$, the {\em restriction} of $\xx$ to $I$ (denoted $\xx^I$) is the vector whose entries are those from $\xx$ taken from positions labelled by $I$.  For example, if $\xx=(1,4,5,2)$ and $I=\{1,3\}$, then $\xx^I = (1,5)$.

%Concatenation
Let $\aa = (a_1,a_2,\ldots,a_m)$ and $\bb = (b_1,b_2,\ldots,b_n)$ be an $m$-tuple and an $n$-tuple respectively. Define the {\em concatenation} of $\aa$
and $\bb$ to be the $(m+n)$-tuple
\[
\cat(\aa,\bb) = (a_1,a_2,\ldots,a_m,b_1,b_2,\ldots,b_n).
\] 
%sum
%weight
The {\em sum} of a vector is simply the sum of its entries, while the {\em weight} of a vector is the number of non-zero entries.  (Note that if a vector's entries are either 0 or 1, its sum and weight are equal.)

Now suppose $\AA=(A_1,A_2,\ldots,A_m)$ and $\BB=(B_1,B_2,\ldots,B_m)$ are $m$-tuples of sets.  We write $\AA\subseteq\BB$ to mean that $A_i\subseteq B_i$ for all $i\in\{1,2,\ldots,m\}$, and say {\em $\AA$ is contained in $\BB$}.  The operations of restriction and concatenation are defined for vectors of sets in the same way they were for vectors of integers.

For any set $X$, we use the notation ${X \choose k}$ to denote the set of all $k$-subsets of $X$.  (Thus if $X$ is finite and has size $n$, then the size of ${X \choose k}$ is ${n \choose k}$.)  If we have an $m$-tuple of sets $\XX=(X_1,X_2,\ldots,X_m)$ and an $m$-tuple of integers $\kk=(k_1,k_2,\ldots,k_m)$, define
\[ {\XX \choose \kk} = {X_1 \choose k_1} \times {X_2 \choose k_2} \times \cdots \times {X_m \choose k_m}. \]
So a member of ${\XX \choose \kk}$ consists of an $m$-tuple of finite sets, of sizes $(k_1,k_2,\ldots,k_m)$.

Other pieces of notation will be defined as and when required.

\subsection{Definition and examples}

Suppose $v,k,t,\lambda$ are integers where $v \geq k \geq t \geq 1$ and $\lambda \geq 1$.  Let $\vv = (v_1,v_2,\ldots,v_m)$ be an $m$-tuple of positive integers with sum $v$, and let $\kk=(k_1,k_2,\ldots,k_m)$ be an $m$-tuple of positive integers with sum $k$, and where $\kk \leq \vv$.

Now let $\XX = (X_1,X_2,\ldots,X_m)$ be an $m$-tuple of pairwise disjoint sets, where $|X_i|=v_i$.  Let $\tt=(t_1,t_2,\ldots,t_m)$ be an $m$-tuple of {\em non-negative} integers.  We say $\tt$ is {\em $(\kk,t)$-admissible} if $\tt \leq \kk$ and $\sum t_i=t$.  In a similar vein, if $\TT=(T_1,T_2,\ldots,T_m)$ is an $m$-tuple of disjoint sets, we say that $\TT$ is {\em $(\vv,\kk,t)$-admissible} if each $T_i$ is a $t_i$-subset of $X_i$, where $(t_1,t_2,\ldots,t_m)$ is $(\kk,t)$-admissible.  (Note that since $t_i$ is allowed to be zero, the corresponding set $T_i$ is allowed to be empty.)

\begin{defn} \label{defn:GC} %%\marginpar{\sf [MODIFIED]}
Suppose $\vv,\kk,t,\lambda,\XX$ are as above.  Then a {\em generalized covering design}, ${\rm GC}_\lambda(\vv,\kk,t)$, is a family $\mathcal{B}$ of elements of ${\XX \choose \kk}$, called {\em blocks}, with the property that every $\TT=(T_1,T_2,\ldots,T_m)$ which is $(\vv,\kk,t)$-admissible is contained in at least $\lambda$ blocks in $\mathcal{B}$.
\end{defn}

We call $X = X_1 \dot\cup X_2 \dot\cup \cdots \dot\cup X_m$ the {\em point set} of the generalized covering design; one can think of $\XX$ as being a partition of the point set $X$.  %%\marginpar{\sf [ADDITION]}
However, by an abuse of notation, we usually label the elements of each $X_i$ as $\{1,2,\ldots,v_i\}$.

We remark that our definition of a generalized covering design is identical to Cameron's definition of a generalized $t$-design, except his definition requires ``exactly $\lambda$'', where ours requires ``at least $\lambda$''.  As with covering arrays, we call the parameter $t$ the {\em strength} of the design.  Now, Cameron's generalized $t$-designs are a common generalization of $t$-designs and orthogonal arrays, whereas our designs are a common generalization of covering designs and covering arrays, as we show below.

\begin{prop} \label{prop:covdes}
Suppose $\vv=(v)$ and $\kk=(k)$.  Then a ${\rm GC}_\lambda(\vv,\kk,t)$ is a $(v,k,t)_\lambda$-covering design.
\end{prop}

\begin{prop} \label{prop:covarr}
Suppose $\vv=s^k$ and $\kk=1^k$.  Then a ${\rm GC}_\lambda(\vv,\kk,t)$ (with $N$ blocks) is equivalent to a covering array $CA_\lambda(N;k,s,t)$.
\end{prop}

\proof Since $\vv=s^k=(s,s,\ldots,s)$, we can suppose $X_1,X_2,\ldots,X_k$ are disjoint copies of an $s$-set $A$, which we regard as our alphabet.  Also, since $\kk=1^k=(1,1,\ldots,1)$, each block consists of exactly one element from each copy of $A$, so we can regard this as a $k$-tuple of elements of $A$.  Put these $k$-tuples as the rows of an array.

Now, each $(\kk,t)$-admissible $k$-tuple $\tt$ will be a 0/1 vector of weight $t$, so represents a $t$-subset of columns of the array.  Thus, for a given $\tt$, the corresponding $(\vv,\kk,t)$-admissible $k$-tuples $T$ will represent all possible combinations of entries in those $t$ columns.  Consequently, each such combination must appear at least once, and so the array is a $CA(N;k,s,t)$. 

Carrying out the reverse of the process described above shows that a covering array $CA(N;k,s,t)$ gives a ${\rm GC}_\lambda(\vv,\kk,t)$ (with $N$ blocks), as required. \endproof

Of course, there are many examples which are neither covering designs nor covering arrays.  The following is a very basic example of such an object.

\begin{example} \label{example:GCDnew}
Let $\vv=(4,2,2)$, $\kk=(2,1,1)$ and $t=2$.  Then the following is a ${\rm GC}(\vv,\kk,2)$:
\[ \begin{array}{lll}
( \{12\}, & \{1\}, & \{1\} ) \\
( \{13\}, & \{1\}, & \{2\} ) \\
( \{14\}, & \{2\}, & \{1\} ) \\
( \{23\}, & \{2\}, & \{2\} ) \\
( \{24\}, & \{1\}, & \{2\} ) \\
( \{34\}, & \{2\}, & \{1\} )
\end{array} \]
The possible admissible vectors $\tt$ are $(2,0,0)$, $(1,1,0)$, $(1,0,1)$ and $(0,1,1)$.  For $\tt=(2,0,0)$, we are required to cover all possible pairs from $\{1,2,3,4\}$ in the first column.  For $\tt=(1,1,0)$, each symbol from $\{1,2,3,4\}$ must appear in the first part of a block with each possible symbol from $\{1,2\}$ in the second part.  The case $\tt=(1,0,1)$ works similarly.  For $\tt=(0,1,1)$ each of the ordered pairs $(1,1)$, $(1,2)$, $(2,1)$, $(2,2)$ must occur in the final two parts of some block.  It is a straightforward exercise to verify that all the possibilities are covered.
\end{example}

As is common in the study of covering problems, from now on we will consider only the case $\lambda=1$, and drop the subscript $\lambda$ from our notation.  Also we note that, by the same argument as Proposition~\ref{prop:covarr}, a mixed covering array ${\rm MCA}(N;k,\vv,t)$ is equivalent to a ${\rm GC}_\lambda(\vv,\kk,t)$ where $\kk=(1,1,\ldots,1)$ (and with $N$ blocks).

As was the case with ordinary covering designs, it is trivial to show that generalized covering designs ${\rm GC}(\vv,\kk,t)$ always exist, provided that $\vv\geq\kk$, simply by taking the collection of all possible blocks and discarding any that are unnecessary.  So the interesting questions, as with all covering problems, are to find bounds on the minimal size of a ${\rm GC}(\vv,\kk,t)$, and find constructions which meet (or come close to) these bounds.  Borrowing notation from the study of (ordinary) covering designs, we have the following.

\begin{notation}
The {\em covering number} $C(\vv,\kk,t)$ denotes the smallest possible size of a generalized covering design ${\rm GC}(\vv,\kk,t)$.
\end{notation}

As with ordinary covering designs, we say $C(\vv,\kk,t)=0$ if no design exists, or if there are no possible $(\vv,\kk,t)$-admissible vectors of sets.

We conclude this section by giving some basic results about generalized covering designs, and the corresponding covering numbers.  Our first result can be used to obtain a recursive bound on covering numbers.

\begin{prop}\label{prop:reducet}
For $t\geq 2$, any ${\rm GC}(\vv, \kk,t)$ is a ${\rm GC}(\vv,\kk,t-1)$, and so 
\[ C(\vv,\kk,t) \geq C(\vv,\kk,t-1). \]
\end{prop}

\proof Let $\TT$ be a $(\vv,\kk,t)$-admissible vector.  Then for any ${\rm GC}(\vv, \kk,t)$ there is a block in the design that contains $\TT$. Clearly, any $m$-tuple $\SS$ with $\SS \subseteq \TT$ is also contained in this block. Since any $(\vv,\kk,t-1)$-admissible vector of sets is contained in some $(\vv,\kk,t)$-admissible vector, the result follows. \endproof

Our next result shows that the case where $t=1$ is very easy.  In this case, all we require is enough blocks so that each symbol from each part appears at least once in the corresponding part of a block.

\begin{prop} \label{prop:tis1}
Where $\vv=(v_1,v_2,\ldots,v_m)$ and $\kk=(k_1,k_2,\ldots,k_m)$, we have
\[
C(\vv, \kk, 1) = \max_{i\in \{1,\ldots,m\}} \left\lceil \frac{v_i}{k_i} \right\rceil.
\]
\end{prop}

\section{The case of strength 2: Clique coverings} \label{section:clique}

A {\em clique covering} of a graph $G$ is a family of complete subgraphs $G_1,\ldots,G_N$ of $G$, called {\em cliques}, with the property that every edge of $G$ appears in at least one of $G_1,\ldots,G_N$.  We remark that the general clique covering problem allows cliques of different sizes.  If we require that all the cliques are {\em $k$-cliques} (i.e.\ they all have exactly $k$ vertices), then we have a {\em $k$-uniform} clique covering.

As with covering designs, the question of the existence of a clique covering of a graph is trivial: the set of all edges forms a 2-uniform clique covering.  On the other hand, for $k>2$, $k$-uniform clique coverings do not exist in general, and their existence depends on the structure of the graph.
	Ideally, we would like a clique covering of $G$ to contain the smallest number of cliques possible; this number is called the {\em clique covering number} of $G$, denoted ${\rm cc}(G)$.  The first result on clique coverings is attributed to Hall \cite{Hall41}, who showed that if a graph has $n$ vertices, then ${\rm cc}(G) \leq \lfloor n^2/4 \rfloor$. 	This was later extended by Erd\H{o}s {\em et al.}\ \cite{Erdos66} who showed the same bound holds when every edge of the graph appears in {\em exactly} one clique of the covering.  Other results on clique coverings can be found in a survey paper by Monson {\it et al.}\ \cite{Monson95}.

	Many techniques and concepts used to construct clique coverings and generate bounds on the clique covering number can be extended to generalized covering designs of strength $2$ and bounds on $C({\bf v},{\bf k},2)$. 	These ideas will be pursued in this section of the paper.  For example, one such concept is that of {\it equivalent vertices} in a graph explored by Gy\'arf\'as \cite{Gyarfas90}, which we apply to generalized covering designs in Section $3.3$.

In the case where $t=2$, both covering designs and covering arrays are equivalent to clique coverings of particular graphs, as we will now explain.

\begin{prop} \label{prop:covdesclique}
A $(v,k,2)$-covering design is a covering of the complete graph $K_v$ by $k$-cliques $G_1,\ldots,G_N$, such that every edge of $K_v$ appears in at least one of $G_1,\ldots,G_N$.
\end{prop}

For instance, a Steiner triple system on $v$ points, which is an example of a $(v,3,2)$-covering design, is a partition of $K_v$ into 3-cliques (so each edge appears exactly once).

In order to explain how covering arrays give clique coverings, we need the following definition.

\begin{defn} \label{defn:multipartite}
Suppose $\mathbf{n}=(n_1,n_2,\ldots,n_m)$ is a sequence of positive integers, and let $V_1,V_2,\ldots,V_m$ be disjoint sets of sizes $n_1,n_2,\ldots,n_m$ respectively.  Then the {\em complete multipartite graph} $K_\mathbf{n}$ has vertex set $V_1\dot\cup V_2 \dot\cup \cdots \dot\cup V_m$, and two vertices $u\in V_i$, $v\in V_j$ are adjacent if and only if $i\neq j$.
\end{defn}

Alternatively, the complete multipartite graph $K_\mathbf{n}$ is the complement of the disjoint union of complete graphs with vertex sets $V_1,V_2,\ldots,V_m$.  In particular: if $\mathbf{n}=(n)$, we have a empty graph; if $\mathbf{n}=(n_1,n_2)$ we have a complete bipartite graph; and if $\mathbf{n}=(1,1,\ldots,1)$ (with $m$ entries) we have the complete graph $K_m$.

\begin{prop} \label{prop:CAclique}
A covering array ${\rm CA}(N;k,s,2)$ is equivalent to a covering of the complete multipartite graph $K_\vv$ by $k$-cliques $K_k$, where $\vv$ has $k$ entries, all equal to $s$.
\end{prop}

\proof Each row of a ${\rm CA}(N;k,s,2)$ defines a clique, where the symbol in position $i$ specifies the vertex chosen from the $i^{\textrm{th}}$ part of the graph $K_\vv$.  In any pair of columns of the array, every pair of symbols must be used at least once; consequently, every edge of $K_\vv$ appears in at least one clique.  On the other hand, given a clique covering of this type, we can construct a covering array by reversing this process.
\endproof

For further details, see Danziger {\em et al.}\ \cite{CAFE}, Maltais \cite{Maltais}, or Ronneseth and Colbourn \cite{Ronneseth09}.  We remark that mixed covering arrays (with $t=2$ and $\lambda=1$) 
%%\marginpar{\sf [MODIFIED]}
can also be interpreted as clique coverings of complete multipartite graphs in exactly the same way, but where the parts may have different sizes.

In order to describe generalized covering designs ${\rm GC}(\vv,\kk,2)$ in terms of clique coverings, we need to define a suitable graph.  To do this, we need the following graph-theoretical idea.

\begin{defn} \label{defn:join}
Let $G_1=(V_1,E_1)$ and $G_2=(V_2,E_2)$ be graphs with $V_1\cap V_2=\emptyset$.  Then the {\em join} of $G_1$ and $G_2$, denoted $G_1+G_2$, is the graph with vertex set $V_1\cup V_2$, and whose edge set is $E_1 \cup E_2 \cup \{xy\, :\, x\in V_1, y\in V_2\}$.
\end{defn}

For example, the join of two complete graphs is also complete, and the join of two empty graphs is a complete bipartite graph.  We note that this can be extended to a join of any number of graphs, and that this operation is associative.

Now suppose that $\vv=(v_1,v_2,\ldots,v_m)$  and $\kk=(k_1,k_2,\ldots,k_m)$ are vectors of positive integers with $\kk\leq\vv$.  Let
\[ H_i=\left\{\begin{array}{cl}
\overline{K_{v_i}}, & \mbox{if $k_i=1$,}\\
K_{v_i}, & \mbox{if $k_i\geq 2$,}\\
\end{array}\right. \]
where $\overline{K_{v_i}}$ represents the complement of $K_{v_i}$ (that is, the empty graph).
Form the graph 
\[ G_{\vv,\kk} = H_1+\cdots+H_m \]
consisting of the join of the graphs $H_i$ such that $G_{\vv,\kk}$ has vertex set $V=\bigcup_i X_i$, where $|X_i|=v_i$ and each $X_i$ is the set of vertices of the corresponding $H_i$.

There are two important special cases.  First, if all $k_i=1$, then $G_{\vv,\kk}$ is precisely the complete multipartite graph $K_\vv$.  Second, if all $k_i\geq 2$, then $G_{\vv,\kk}$ is isomorphic to the complete graph $K_v$ (where $v=v_1+v_2+\cdots+v_m$), but where the vertex set has a specific partition into parts of sizes $v_1,v_2,\ldots,v_m$.

\begin{thm} \label{thm:GCgraph2}
Let $G_{\vv,\kk}$ be the graph described above.  Then a generalized covering design ${\rm GC}(\vv,\kk,2)$ is equivalent to an edge covering of $G_{\vv,\kk}$ using a collection of cliques of size $k$, with the property that for each clique in the covering, $k_i$ vertices of the clique come from the set $X_i$ (for each $i$).
\end{thm}

\proof Suppose $\mathcal{D}$ is a ${\rm GC}(\vv,\kk,2)$.  Now, from each block in $\mathcal{D}$, we can easily construct a clique in $G_{\vv,\kk}$.  Think of the vertex set of $G_{\vv,\kk}$ as being the point set of $\mathcal{D}$, namely $X_1\dot\cup X_2 \dot\cup \cdots \dot\cup X_m$.  A block of $\mathcal{D}$ contains $k_i$ points from each part $X_i$, and the subgraph of $G_{\vv,\kk}$ induced by the corresponding vertices is necessarily a clique.  Now, the admissible vectors $\tt$ have two possible forms: (i) a vector with two entries~1 and the rest~0; and (ii) a vector with a single entry of 2, and all other entries~0.  Type (i) vectors ensure that this collection of cliques covers all edges between parts.  Type (ii) vectors are only possible with the single~2 in position~$i$ with $k_i\geq 2$; these vectors ensure that these cliques cover all edges within each part with $k_i\geq 2$.  Consequently, these cliques cover all the edges of $G_{\vv,\kk}$, and so we have a clique covering.

A similar argument works in the reverse direction: given a clique covering of $G_{\vv,\kk}$ with the specified form, each clique gives us a block, and the way $G_{\vv,\kk}$ was constructed ensures these blocks form a ${\rm GC}(\vv,\kk,2)$.  \endproof

\begin{example} \label{example:cliqueGC}
Recall Example~\ref{example:GCDnew}, where we saw a ${\rm GC}(\vv,\kk,2)$ with $\vv=(4,2,2)$ and $\kk=(2,1,1)$.  By Theorem~\ref{thm:GCgraph2}, this can be viewed as a clique covering of the appropriate graph $G_{\vv,\kk}$ using copies of $K_4$, as shown in Figure~\ref{figure:cliqueGC}.
\setlength{\unitlength}{1cm}
\begin{figure}[hbt]
\centering
\subfigure[The graph $G_{\vv,\kk}$, where $\vv=(4,2,2)$ and $\kk=(2,1,1)$.]{%
\begin{picture}(15,3.5)
\put(5.5,0.5){\includegraphics[width=4cm]{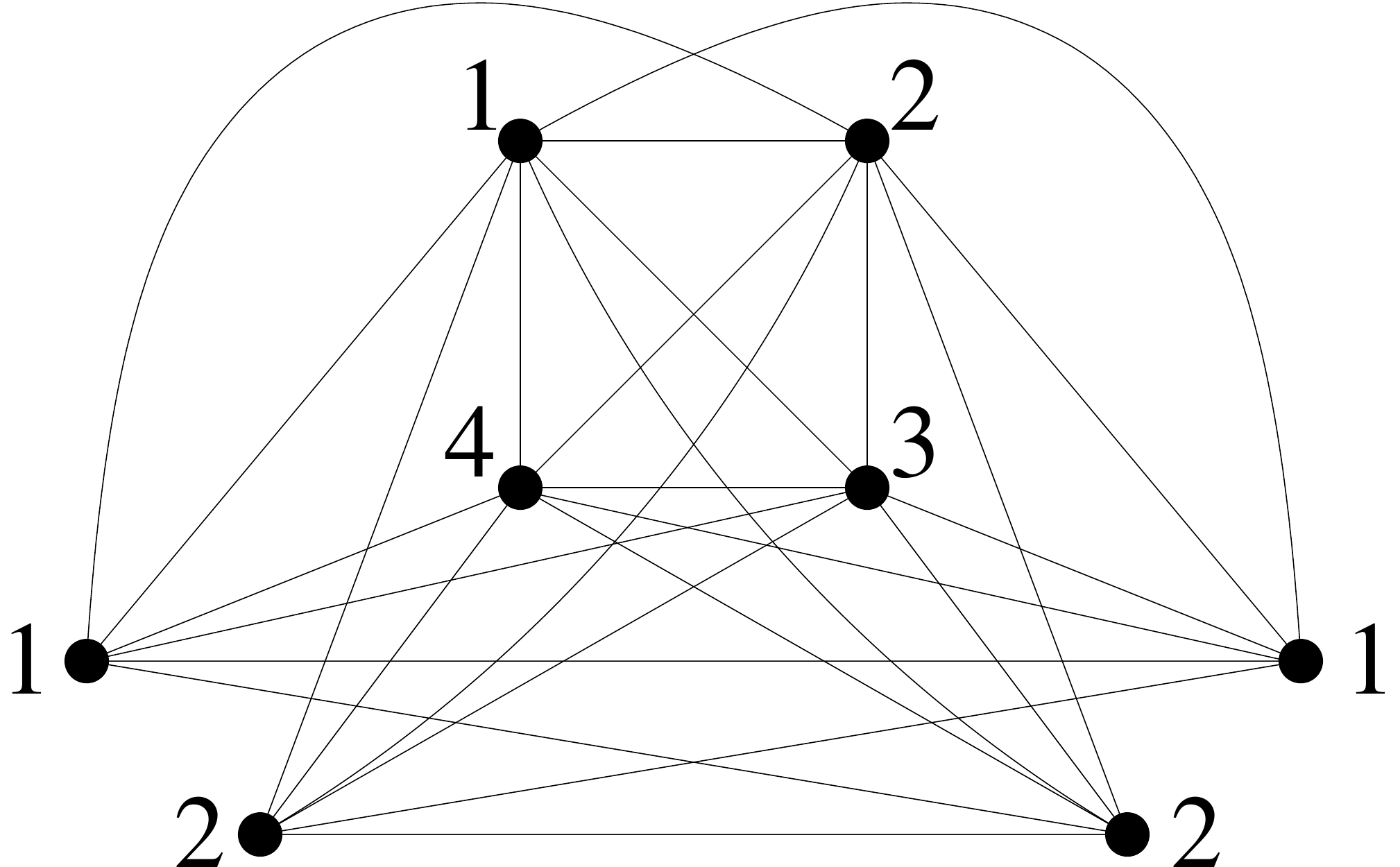}}
\end{picture}
}
\subfigure[A clique covering of $G_{\vv,\kk}$.]{%
\begin{picture}(14.5,8)
\put(0,5){\includegraphics[width=4cm]{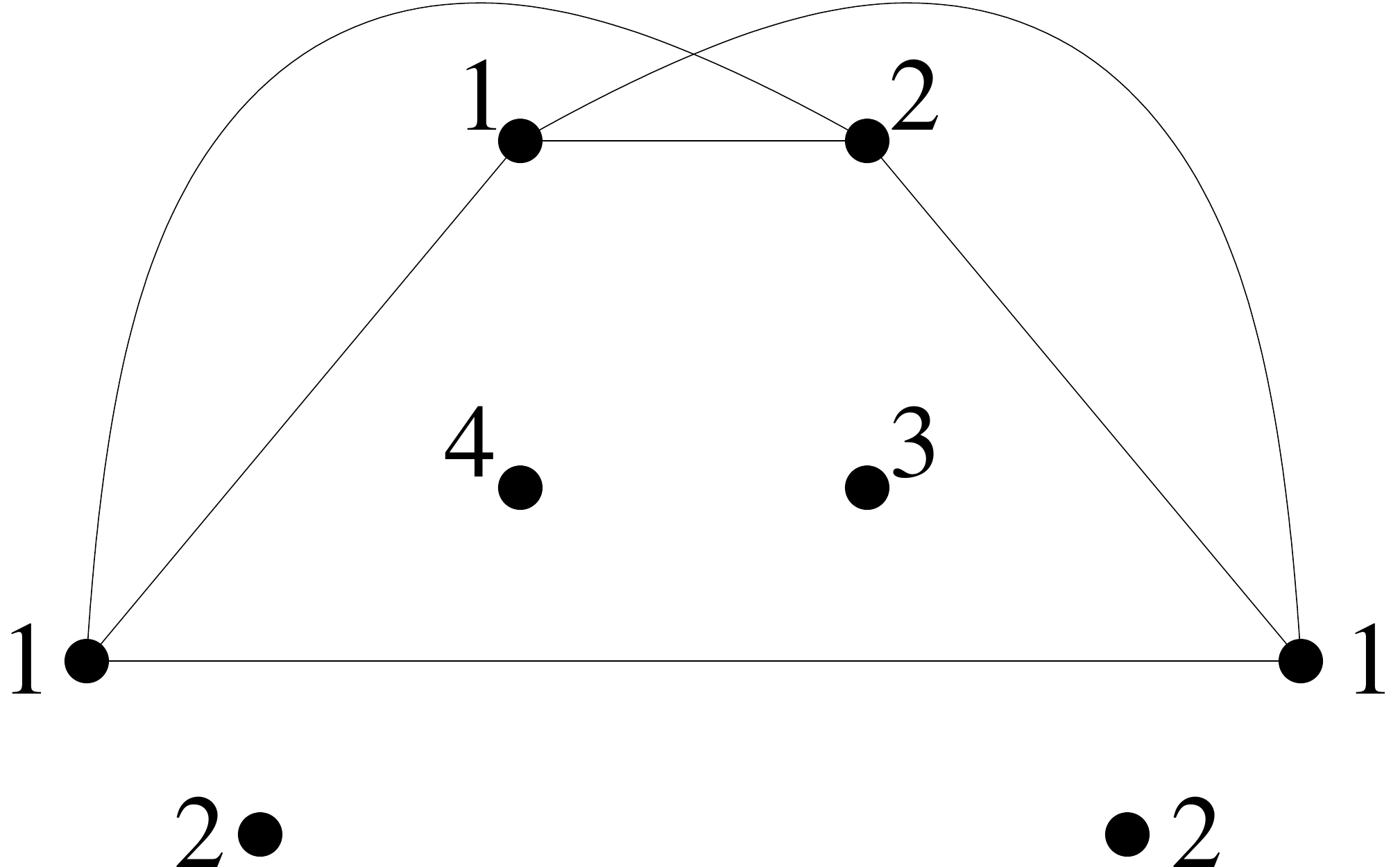}}
\put(5,5){\includegraphics[width=4cm]{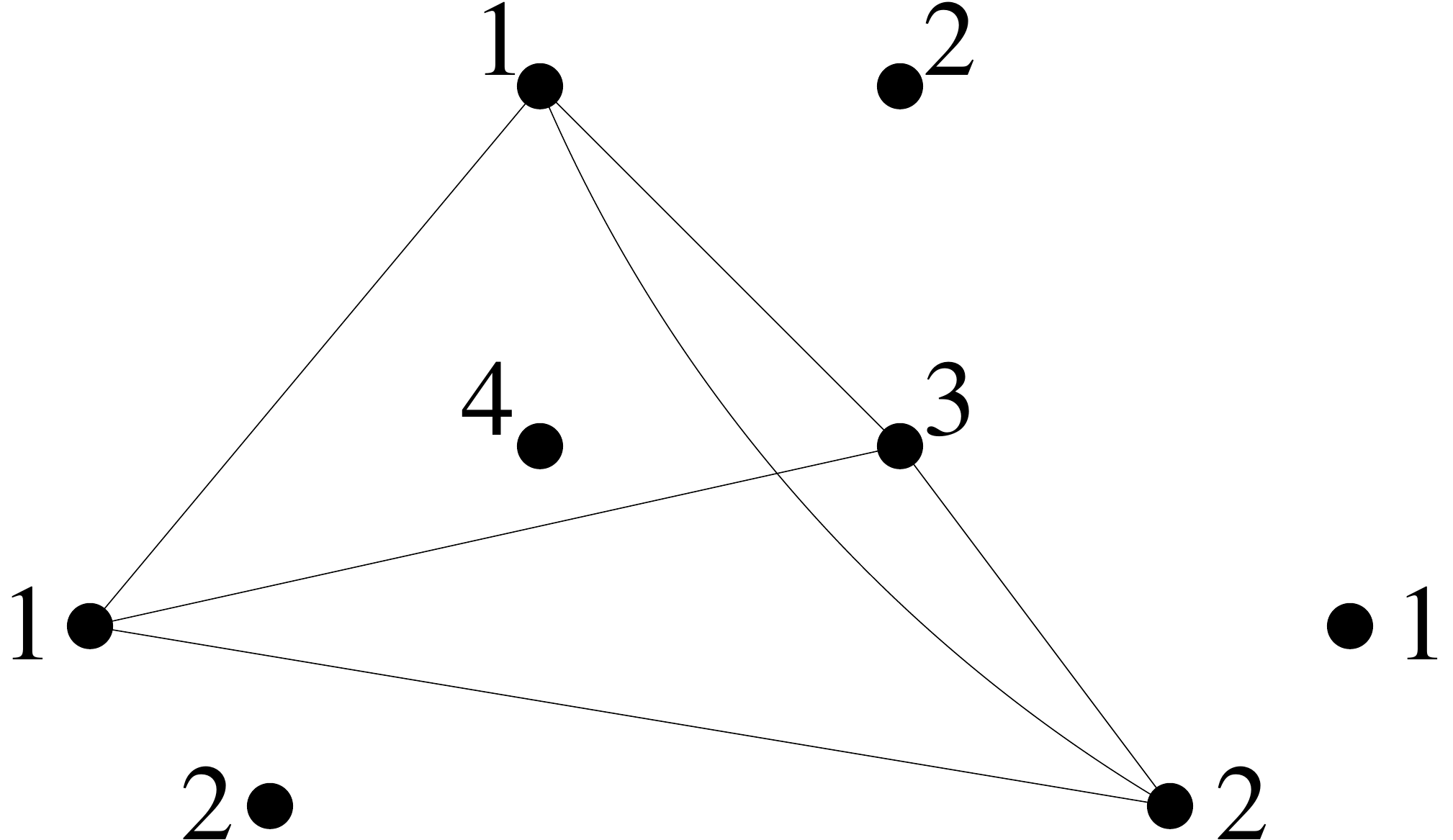}}
\put(10,5){\includegraphics[width=4cm]{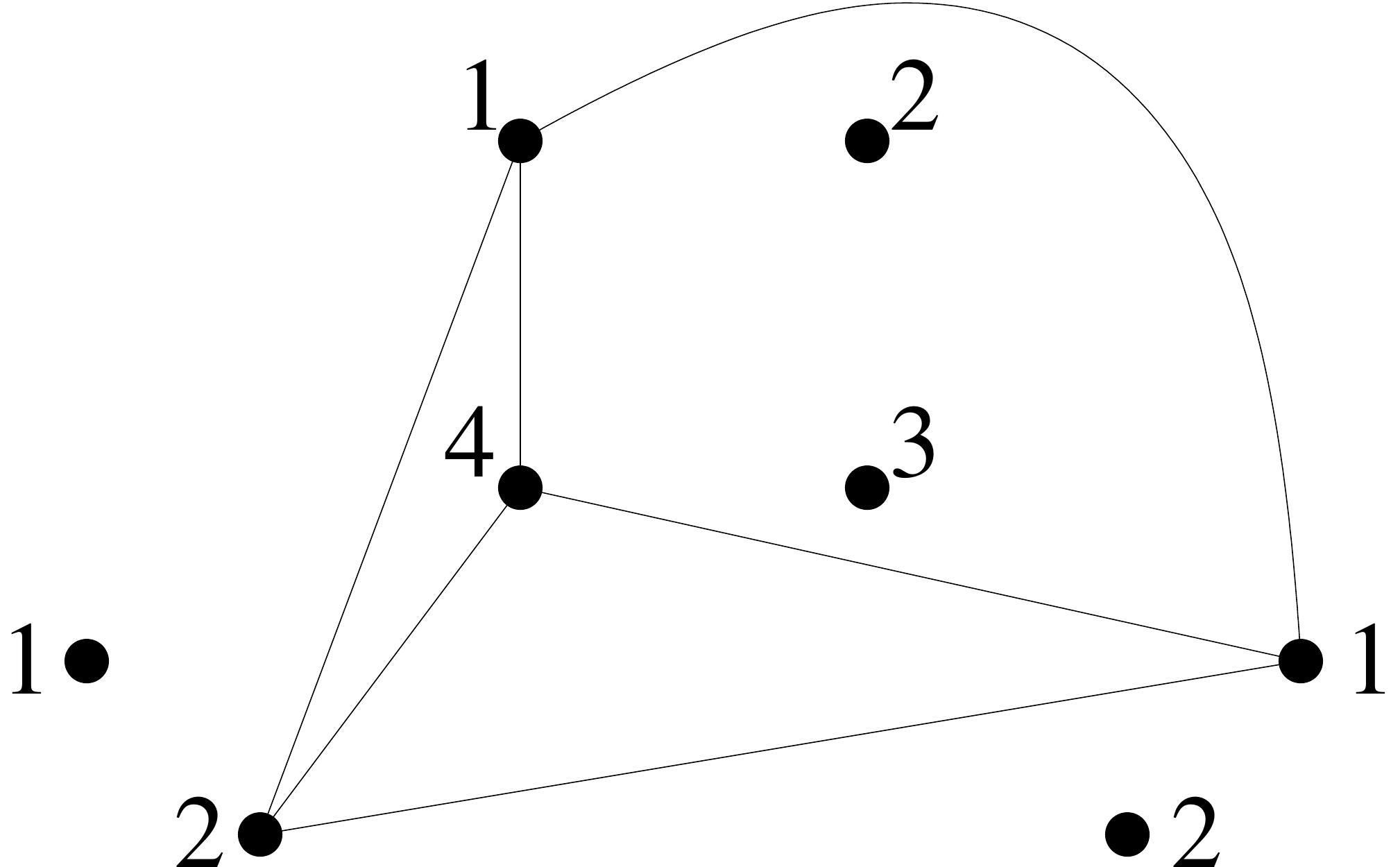}}
\put(0,1){\includegraphics[width=4cm]{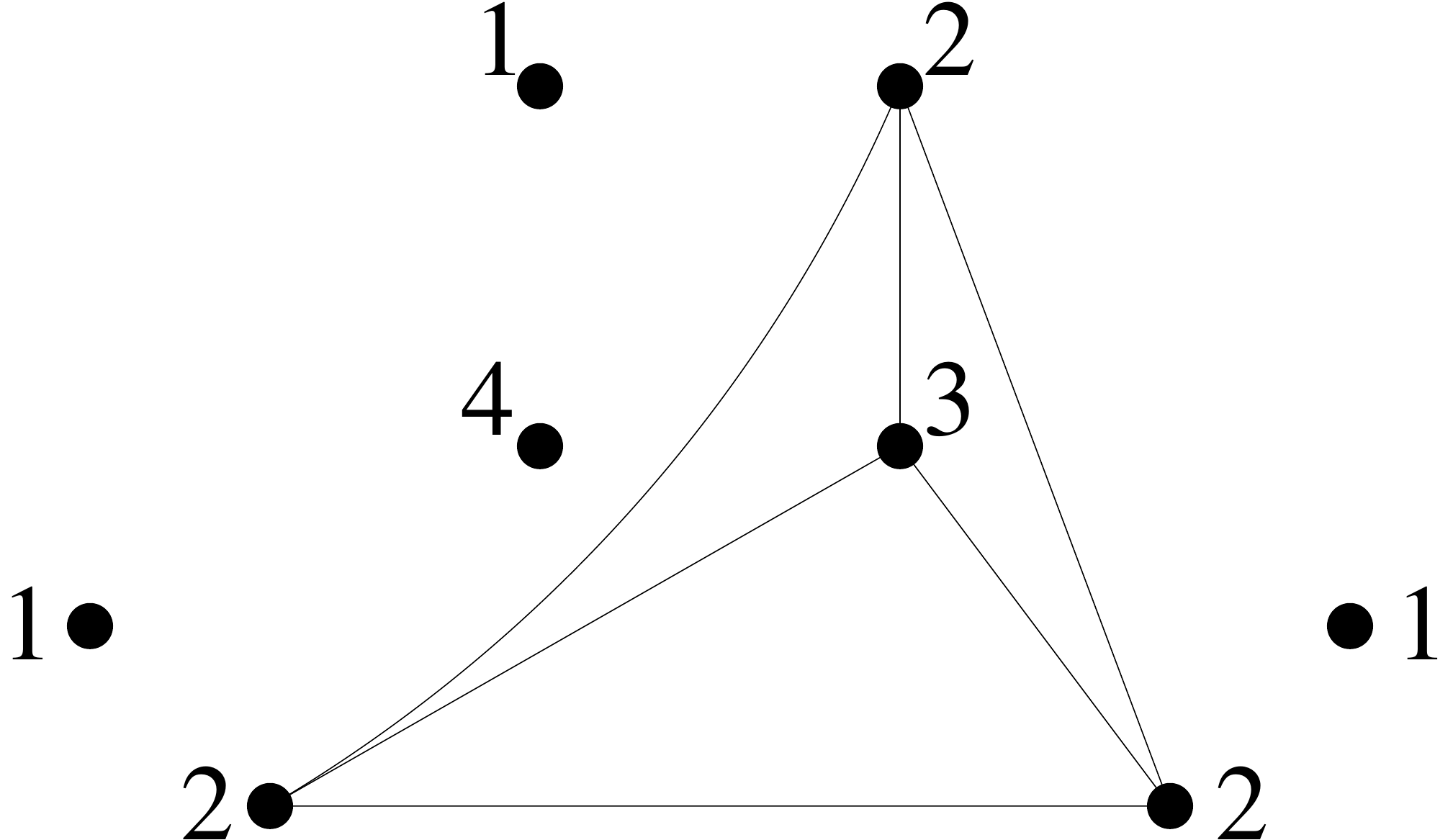}}
\put(5,1){\includegraphics[width=4cm]{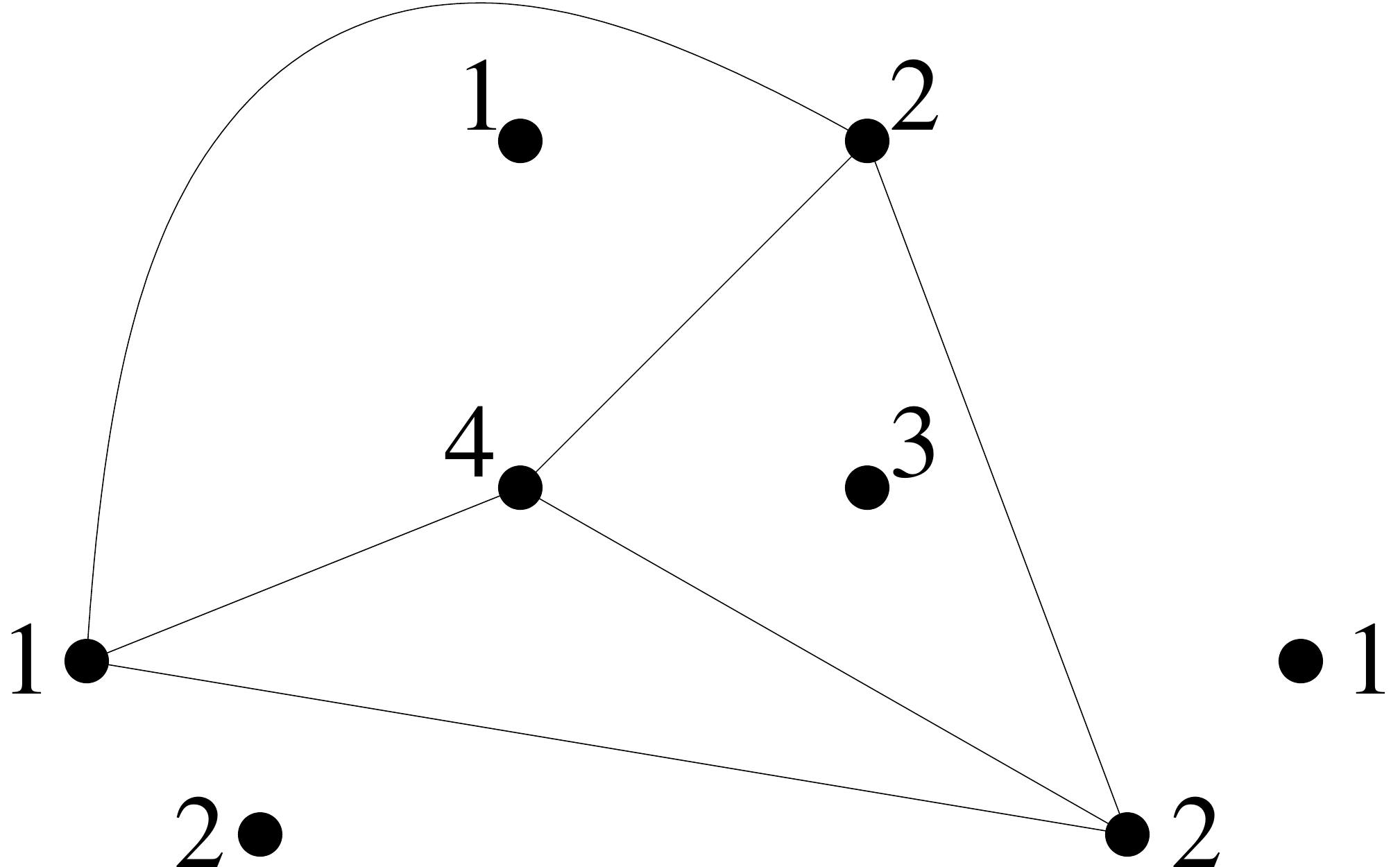}}
\put(10,1){\includegraphics[width=4cm]{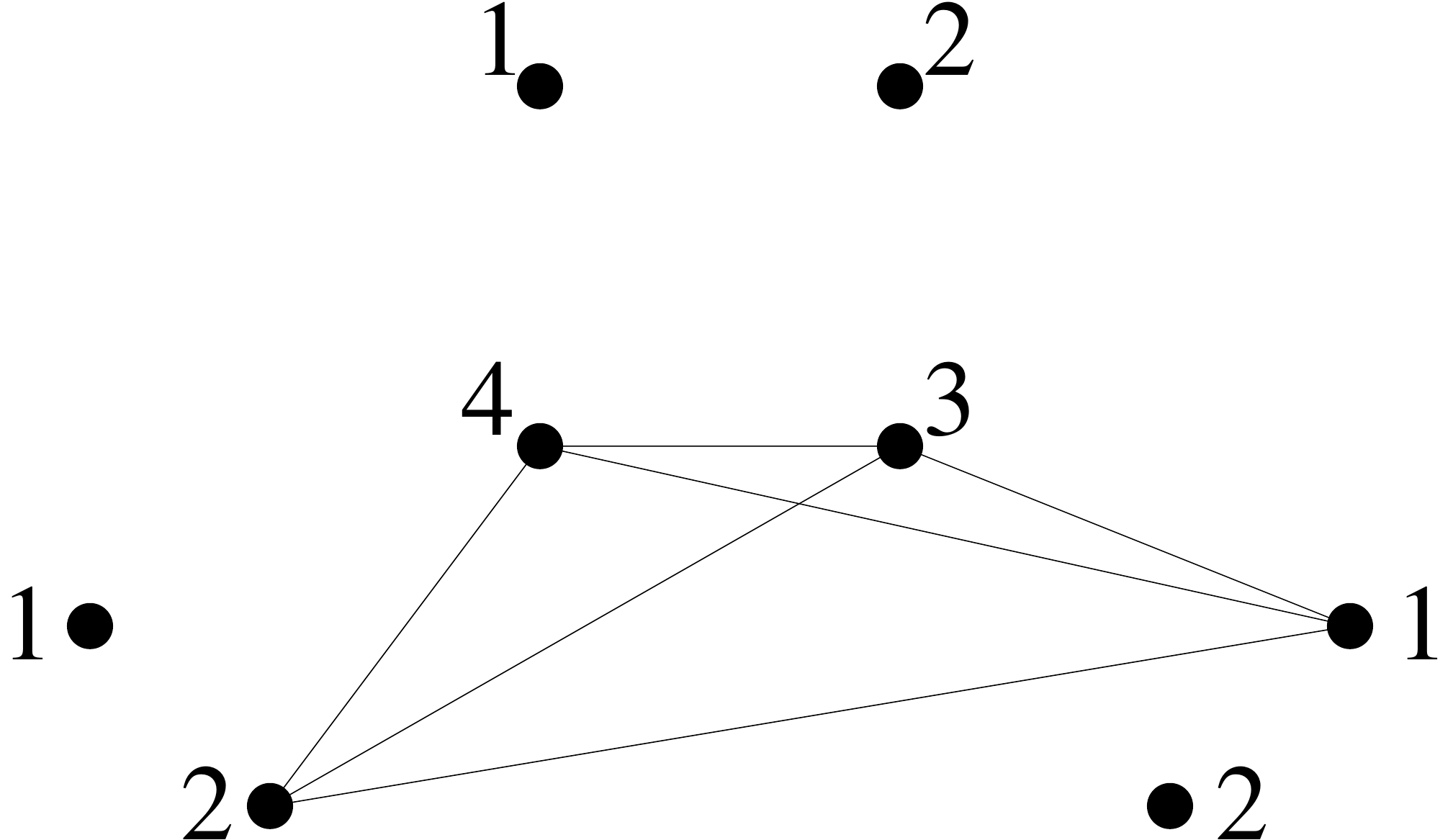}}

\put(0.8,4.5){\footnotesize $( \{12\}, \; \{1\}, \; \{1\} )$}
\put(5.8,4.5){\footnotesize $( \{13\}, \; \{1\}, \; \{2\} )$}
\put(10.8,4.5){\footnotesize $( \{14\}, \; \{2\}, \; \{1\} ) $}
\put(0.8,0.5){\footnotesize $( \{23\}, \; \{2\}, \; \{2\} )$}
\put(5.8,0.5){\footnotesize $( \{24\}, \; \{1\}, \; \{2\} )$}
\put(10.8,0.5){\footnotesize $( \{34\}, \; \{2\}, \; \{1\} )$}
\end{picture}
}
\caption{The generalized covering design from Example~\ref{example:GCDnew}, viewed as a clique covering. \label{figure:cliqueGC} }
\end{figure}
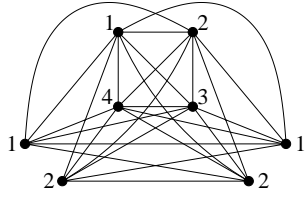
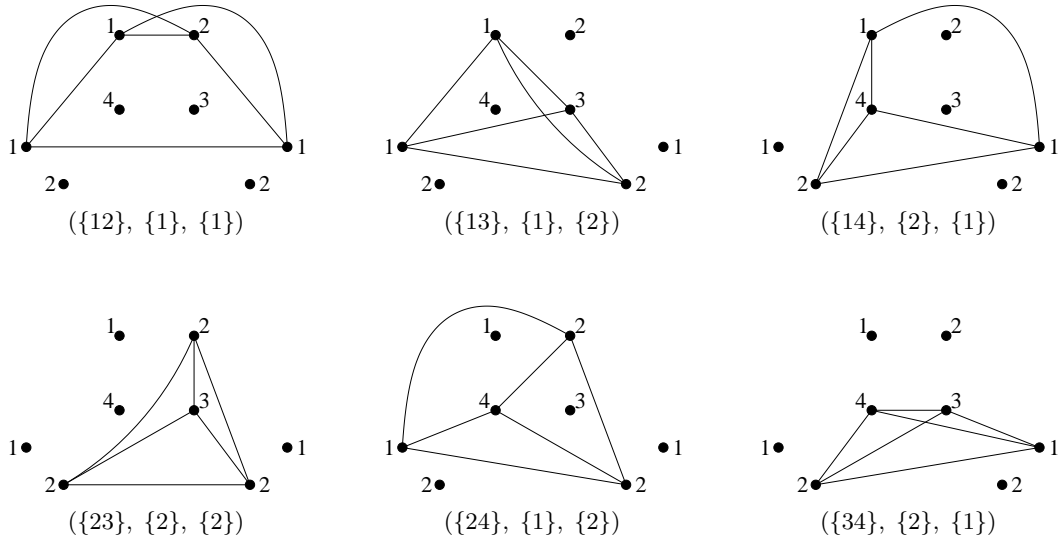
\end{example}

We can use this interpretation to obtain various bounds on covering numbers $C(\vv,\kk,2)$, as well as to obtain methods of constructing generalized covering designs ${\rm GC}(\vv,\kk,2)$.  In some cases, we can use these results to obtain designs which are optimal.

\subsection{Counting edges}
The simplest bound on $C(\vv,\kk,2)$ can be obtained by counting edges.  When covering a graph $G$ (with $|E(G)|$ edges) by a collection of subgraphs isomorphic to $H$, each of which has $|E(H)|$ edges, then clearly the minimum number of subgraphs required is $\lceil |E(G)|/|E(H)| \rceil$.  This observation then provides the following bound:
\[ C(\vv,\kk,2)\geq \left\lceil  \frac{|E(G_{\vv,\kk})|}{|E(K_k)|}  \right\rceil.\]
Using the interpretation in Theorem~\ref{thm:GCgraph2}, the number of edges in $G_{\vv,\kk}$ is
\[ |E(G_{\vv,\kk})|=\binom{v}{2}-\sum_{k_i=1}\binom{v_i}{2}, \]
where the summation is over the indices $i$ such that $k_i=1$, and we use the convention that $\binom{1}{2}=0$.  However, an alternative way to write $|E(G_{\vv,\kk})|$ is
\[ |E(G_{\vv,\kk})|=\sum_{k_i\neq 1} \binom{v_i}{2}+\sum_{i\neq j}v_i v_j, \]
as the first summation counts the edges in each $H_i$ and the second summation counts the edges joining each $H_i$ to each $H_j$ (for $i\neq j$).  Now, since the number of edges in a clique of size $k$ is $|E(K_k)|=\binom{k}{2}$, we have shown the following.

\begin{prop} \label{prop:naive}
Using the notation above, 
\[ C(\vv,\kk,2)\geq 
\left\lceil \frac{ \binom{v}{2}-\sum_{k_i=1}\binom{v_i}{2} }{ {k \choose 2} } \right\rceil = 
\left\lceil \frac{ \sum_{k_i\neq 1} \binom{v_i}{2}+\sum_{i\neq j}v_i v_j }{ {k \choose 2} } \right\rceil. \]
\end{prop}

\subsection{Restriction}
%In this subsection, we give an operation which helps construct generalized covering designs, and which works as follows.  Suppose we have a vector $\xx=(x_1,\ldots,x_m)$, and let $I$ be a non-empty subset of the positions $\{1,\ldots,m\}$.  Then the {\em restriction} of $\xx$ to $I$, denoted $\xx^I$, is the vector with $|I|$ entries, consisting of those entries of $\xx$ labelled by the elements of $I$.
Recall the notion of {\em restriction} which we introduced earlier.  In this subsection, we will show how this operation can be used to construct generalized covering designs, and to obtain bounds on $C(\vv,\kk,2)$.

\begin{prop} \label{prop:restriction}  %%\marginpar{\sf [MODIFIED]}
Let $\mathcal{D}$ be a generalized covering design ${\rm GC}(\vv,\kk,2)$ with $N$ blocks, and suppose $I\subseteq\{1,\ldots,m\}$ is non-empty.  Then, provided $\kk^I\neq(1)$, there exists a ${\rm GC}(\vv^I,\kk^I,2)$ with $N$ blocks, and in particular $C(\vv,\kk,2)\geq C(\vv^I,\kk^I,2)$.
\end{prop}

\proof %%\marginpar{\sf [MODIFIED]}
We exclude the case $\kk^I=(1)$ as no $(v_i,1,2)$-covering design can exist; however, since we adopted the convention that $C(v_i,1,2)=0$, the bound still holds, albeit for trivial reasons.

So we suppose $\kk\neq(1)$ and 
construct the graph $G_{\vv,\kk}$ as in Theorem~\ref{thm:GCgraph2} above.  Clearly, the graph $G_{\vv^I,\kk^I}$ is an induced subgraph of $G_{\vv,\kk}$, obtained simply by removing the parts not indexed by $I$.  Now, each block of $\mathcal{D}$ is a clique in $G_{\vv,\kk}$ with $k_i$ vertices in part $X_i$, and if we restrict these cliques to those parts in $G_{\vv^I,\kk^I}$, every edge of $G_{\vv^I,\kk^I}$ must be covered by at least one clique.  Consequently, these $N$ cliques form a clique covering of $G_{\vv^I,\kk^I}$ with the required property, and thus correspond to a ${\rm GC}(\vv^I,\kk^I,2)$ with $N$ blocks.  The bound $C(\vv,\kk,2)\geq C(\vv^I,\kk^I,2)$ follows immediately.  \endproof

We denote the ${\rm GC}(\vv^I,\kk^I,2)$ obtained from $\mathcal{D}$ in the proof of Proposition~\ref{prop:restriction} by $\mathcal{D}^I$, and call this the {\em restriction of $\mathcal{D}$ to $I$}.  Now, by restricting in all possible ways, we obtain the following lower bound.

\begin{cor} \label{cor:restriction}
Let $\mathscr{I}$ denote the collection of all non-empty subsets of $\{1,\ldots,m\}$.  Then we have
\[ C(\vv,\kk,2)\geq  \max_{I\in \mathscr{I}} C(\vv^I,\kk^I,2). \]
\end{cor}

Of course, the na\"{\i}ve edge-counting bound in Proposition \ref{prop:naive} can then be applied to the right-hand side to obtain a lower bound in terms of the entries of $\vv$ and $\kk$.

Unfortunately, the bound in Corollary \ref{cor:restriction} involves checking an exponential number of cases, so a more practical bound can be obtained just by looking at subsets of size $1$.  In this situation, the restricted design has only one part, and is thus an ordinary $(v_i,k_i,2)$-covering design.  Consequently, we have the following corollary.

\begin{cor} \label{cor:restriction2}
For vectors $\vv=(v_1,\ldots,v_m)$ and $\kk=(k_1,\ldots,k_m)$, we have 
\[ C(\vv,\kk,2) \geq \max_{\substack{1\leq i \leq m \\ k_i\neq 1}} \{ C(v_i,k_i,2) \}. \]
\end{cor}

We were able to exclude indices where $k_i=1$ from the above, as $C(v_i,1,2)=0$.  Of course, if $k_i \geq 2$ for some $i$ (i.e.\ our design is not a (mixed) covering array), we can then use the Sch\"onheim bound (Theorem~\ref{thm:Schonheim}) to obtain
\[ C(\vv,\kk,2) \geq \max_{\substack{1\leq i \leq m \\ k_i\neq 1}} \left\{ \left\lceil \frac{v_i}{k_i} \left\lceil \frac{v_i-1}{k_i-1} \right\rceil \right\rceil \right\}. \]

We also have another straightforward application of restriction, which allows us to safely ignore parts where $v_i=k_i$.

\begin{prop} \label{prop:allvertexparts}
Suppose $\vv=(v_1,\ldots,v_m)$ and $\kk=(k_1,\ldots,k_m)$ are vectors of positive integers with $\vv\geq\kk$.  Let $I$ be the set of indices where $v_i\neq k_i$.  Then, provided that $\kk^I \neq (1)$, we have $C(\vv,\kk,2)=C(\vv^I,\kk^I,2)$.
\end{prop}

\proof We require that $\kk^I \neq (1)$ in order to ensure that the graph $G_{\vv^I,\kk^I}$ is not an empty graph (in which case, we would have $C(\vv^I,\kk^I,2)=C(v_i,1,2)=0$ and the result would not hold).

From Proposition~\ref{prop:restriction}, we know that $C(\vv,\kk,2)\geq C(\vv^I,\kk^I,2)$; to show that equality holds, we will show that given a ${\rm GC}(\vv^I,\kk^I,2)$, we can obtain a ${\rm GC}(\vv,\kk,2)$ with the same number of blocks.  

Suppose without loss of generality that $I=\{1,\ldots,r\}$, and let $\mathcal{E}$ be a ${\rm GC}(\vv^I,\kk^I,2)$, which we treat as a clique covering of $G_{\vv^I,\kk^I}$.  Now, to each clique in $\mathcal{E}$, we add $m-r$ new parts of sizes $v_{r+1},\ldots,v_m$.  In doing so we obtain a ${\rm GC}(\vv,\kk,2)$: (i) every edge we need to cover in parts $1,\ldots,r$ is already covered in $\mathcal{E}$; (ii) any edge between parts $i,j\in\{r+1,\ldots,m\}$ appears in every block, as does any edge within a part $i\in\{r+1,\ldots,m\}$ (where $k_i\geq 1$); (iii) every edge joining a part $i\in\{r+1,\ldots,m\}$ to a vertex $v\in V(G_{\vv^I,\kk^I})$ is covered, as the vertex $v$ must appear in some block of $\mathcal{E}$. \endproof

We remark that Proposition~\ref{prop:allvertexparts} works for arbitrary values of $t\geq 2$, although the argument using clique coverings can't be applied.

\subsection{Equivalence}
The notion of {\em equivalence} is a very useful one, when it comes to both bounds and constructions.  It is inspired by the notion of {\em equivalent vertices} for clique coverings, as studied by Gy\'arf\'as \cite{Gyarfas90}.

As usual, we have vectors $\vv=(v_1,\ldots,v_m)$ and $\kk=(k_1,\ldots,k_m)$ with $\vv\geq\kk$.  We define an equivalence relation $\sim$ on the set of indices $\{1,\ldots,m\}$ by $i\sim j$ if and only if $v_i=v_j$ and $k_i=k_j$.  (Sometimes, it is convenient to talk about $(v_i,k_i)$ being equivalent to $(v_j,k_j)$.)  In order to demonstrate the usefulness of this idea, we have the following lemma.

\begin{lemma} \label{lemma:equiv}
Where $\vv$ and $\kk$ are as above, let $\mathcal{D}$ be a ${\rm GC}(\vv,\kk,2)$ of size $N$.  Define vectors $\vv^\ast=(v_1,\ldots,v_m,w)$ and  $\kk^\ast=(k_1,\ldots,k_m,\ell)$ where $(w,\ell)$ is equivalent to some fixed $(v_i,k_i)$ where $\ell=k_i\geq 2$.  Then there exists a ${\rm GC}(\vv^\ast,\kk^\ast,2)$ of size $N$.
\end{lemma}

\proof We think of $\mathcal{D}$ in terms of a clique covering of $G_{\vv,\kk}$ (cf. Theorem~\ref{thm:GCgraph2}) with vertex set $X_1\cup X_2 \cup \cdots \cup X_m$.  Since $\ell=k_i\geq 2$, we observe that $G_{\vv^\ast,\kk^\ast}$ is precisely the graph obtained from $G_{\vv,\kk}$ by joining a complete graph $K_w$, i.e.\ $G_{\vv^\ast,\kk^\ast} = G_{\vv,\kk}+K_w$, and where $G_{\vv^\ast,\kk^\ast}$ has vertex set $X_1\cup X_2 \cup \cdots \cup X_m \cup Y$ (where $|Y|=w$).  Since $(w,\ell)\sim(v_i,k_i)$, we have a one-to-one correspondence between the vertices in $X_i$ and in $Y$; label these as $X_i=\{x_1,\ldots,x_{v_i}\}$ and $Y=\{y_1,\ldots,y_w\}$ (of course, $w=v_i$).

Choose some $K\in\mathcal{D}$ (i.e.\ a clique in $G_{\vv,\kk}$): without loss of generality, assume this clique contains vertices $\{x_1,\ldots,x_{k_i}\}$ from $X_i$.  Now form a new clique $K^\ast$ by taking all vertices of $K$, together with $\{y_1,\ldots,y_{k_i}\}$.  We claim that the set $\mathcal{D}^\ast = \{ K^\ast \mid K\in \mathcal{D} \}$ is a clique covering of $G_{\vv^\ast,\kk^\ast}$ corresponding to a ${\rm GC}(\vv^\ast,\kk^\ast,2)$.  

To show this, we need to consider any pair of adjacent vertices of $G_{\vv^\ast,\kk^\ast}$, and show that the edge joining them is covered by at least one member of $\mathcal{D}^\ast$.  Now, any edge which lies in the induced subgraph $G_{\vv,\kk}$ is automatically covered (because we had a clique covering of $G_{\vv,\kk}$ to begin with).  Also, any edge $zy_a$ (where $z$ is a vertex of $G_{\vv,\kk}$ not in $X_i$) is covered, as the corresponding edge $zx_a$ is in some clique in $\mathcal{D}$.  Finally, since $k_i\geq 2$, we know that each edge $x_a x_b$ is covered by some clique $K\in\mathcal{D}$.  Thus both the vertices $y_a$ and $y_b$ are in $K^\ast\in\mathcal{D}^\ast$, and so it follows that $K^\ast$ covers the edges $y_a x_a$, $y_a x_b$ and $y_a y_b$.

Thus all types of edge in $G_{\vv^\ast,\kk^\ast}$ are covered by at least one clique in $\mathcal{D}^\ast$, and we are done.  \endproof

The proof of Lemma~\ref{lemma:equiv} demonstrates the reason for the requirement that $k_i\geq 2$ in the statement of the lemma: if $k_i=\ell=1$, we would not be able to cover edges of the form $y_a x_b$ using this construction.  With that in mind, without loss of generality we assume that $\kk=\cat(1^r,\bl)$, where all entries of $\bl$ are at least~2.  Now write $\vv=\cat(\uu,\ww)$, where $u_i=v_i$ for $i=1,\ldots,r$, and where $w_i=v_{i+r}$ for $i=1,\ldots,m-r$.  Then we have the following theorem.

\begin{thm} \label{thm:equiv}
Suppose we have vectors of positive integers $\uu,\ww,\bl$ as above, where $\ell_i\geq 2$ for all $i$.  Let $R$ be a set of equivalence class representatives for $\ww$ and $\bl$ under the relation~$\sim$.  Then we have
\[ C(\cat(\uu,\ww),\cat(1^r,\bl),2) = C(\cat(\uu,\ww^R),\cat(1^r,\bl^R),2). \]
\end{thm}

\proof The inequality $C(\cat(\uu,\ww),\cat(1^r,\bl),2) \geq C(\cat(\uu,\ww^R),\cat(1^r,\bl^R),2)$ follows from Proposition~\ref{prop:restriction}, by restricting to $R$.  The reverse inequality follows by repeatedly applying Lemma~\ref{lemma:equiv} to a minimal ${\rm GC}(\cat(\uu,\ww^R),\cat(1^r,\bl^R),2)$, to obtain a ${\rm GC}(\cat(\uu,\ww),\cat(1^r,\bl),2)$ of size $C(\cat(\uu,\ww^R),\cat(1^r,\bl^R),2)$. \endproof

In the special case where all entries of $\kk$ are at least 2, we have the following straightforward corollary.

\begin{cor} \label{cor:equiv}
Suppose $k_i\geq 2$ for all $i$, and that $R$ is a set of equivalence class representatives for $\vv$ and $\kk$ under the relation $\sim$.  Then we have
\[ C(\vv,\kk,2) = C(\vv^R,\kk^R,2). \]
\end{cor}

\subsection{Point deletion and block expansion}
The operations of point deletion and block expansion give us another way to construct new generalized covering designs from existing ones.  Suppose that $\vv$ and $\kk$ are our usual vectors of integers, where $\vv=(v_1,\ldots,v_m)$ lists the sizes of the sets $\XX=(X_1,\ldots,X_m)$ respectively.  Now, by {\em point deletion} we mean the operation of removing some points from each part $X_i$, to obtain $\hat{\XX} = (\hat{X_1},\ldots,\hat{X_m})$, with sizes $\hat\vv = (\hat{v_1},\ldots,\hat{v_m})$ respectively (so that $\hat\vv \leq \vv$).  A similar operation is {\em block expansion}, where we keep $\vv$ fixed but increase the sizes of our blocks from $\kk$ to $\hat\kk$, where $\kk\leq\hat\kk$.

In the following discussion, we once again think of our ${\rm GC}(\vv,\kk,2)$ in terms of a clique covering of $G_{\vv,\kk}$ (cf. Theorem~\ref{thm:GCgraph2}) with vertex set $X_1\dot\cup X_2 \dot\cup \cdots \dot\cup X_m$.  Our first lemma concerns point deletion.

\begin{lemma} \label{lemma:deletion}
Let $\mathcal{D}$ be a ${\rm GC}(\vv,\kk,2)$ of size $N$, and suppose $\hat\vv$ satisfies $\kk \leq \hat\vv \leq\vv$.  Then there exists a ${\rm GC}(\hat\vv,\kk,2)$ of size at most $N$.
\end{lemma}

\proof Suppose we have deleted the vertex $x\in X_i$.  For each clique in $\mathcal{D}$ which includes the vertex $x$, we remove $x$ and replace it with another vertex $y\in \hat{X_i}$, which is not already in that clique; since $\hat{v_i}\geq k_i$, we know that such a vertex must exist.  We repeat this procedure until we are left only with the vertices in $\hat\XX$.  It is straightforward to see that the cliques obtained will cover all edges in $G_{\hat\vv,\kk}$.

Note that it is possible that this procedure may introduce some repeated blocks, so we remove any duplicates. Thus we can only be sure our ${\rm GC}(\hat\vv,\kk,2)$ has an upper bound of at most $N$ blocks, rather than exactly $N$. \endproof

A counterpart to Lemma~\ref{lemma:deletion} is the following, which is concerned with block expansion.

\begin{lemma} \label{lemma:expansion}
Let $\mathcal{D}$ be a ${\rm GC}(\vv,\kk,2)$ of size $N$, where $k_i\geq 2$ for all $i$, and suppose $\hat\kk$ satisfies $\kk \leq \hat\kk \leq\vv$.  Then there exists a ${\rm GC}(\vv,\hat\kk,2)$ of size at most $N$.
\end{lemma}

\proof Since $k_i\geq 2$ for all $i$, it follows that $G_{\vv,\hat{\kk}}=G_{\vv,\kk}$.  Thus every edge of $G_{\vv,\hat{\kk}}$ is already covered by the cliques in $\mathcal{D}$, which each have $k_i$ vertices chosen from $X_i$ (for every $i$).  Adding extra vertices to each clique, in parts where $\hat{k_i}>k_i$ (so that there are now $\hat{k_i}$ vertices chosen from each $X_i$), does not affect this.

Again, we note that it is possible that this procedure may introduce some repeated blocks, so we remove any duplicates and thus only have an upper bound of $N$ rather than equality.
\endproof

Of course, in both Lemma~\ref{lemma:deletion} and Lemma~\ref{lemma:expansion}, if the generalized design $\mathcal{D}$ that we begin with is optimal, there is no guarantee that the resulting design will be optimal.  However, the two constructions do yield the following important bound.

\begin{thm} \label{thm:delexp}
Suppose we have vectors of positive integers $\vv,\ww,\kk,\bl$, each of length $m$, where $k_i\geq 2$ for all $i$, and which satisfy $\kk \leq \ww \leq \vv$ and $\kk \leq \bl \leq \vv$.  Then we have
\[ C(\vv,\kk,2) \geq C(\ww,\kk,2) \]
and
\[ C(\vv,\kk,2) \geq C(\vv,\bl,2). \]
\end{thm}

We can combine the operations of point deletion and block expansion, and the notion of equivalence, to obtain our next bound, which is perhaps the most useful so far.  Suppose we have our usual vectors $\vv=(v_1,\ldots,v_m)$ and $\kk=(k_1,\ldots,k_m)$.  Define $\vmax = \max v_i$, and let $\vvmax$ denote the vector of $m$ entries, all equal to $\vmax$; in a similar fashion, define $\kkmin$ to be the vector with $m$ entries, all equal to $\kmin = \min k_i$.  Clearly, $\vvmax \geq \vv$ and $\kkmin \leq \kk$.

\begin{thm} \label{thm:minimax}
Suppose we have vectors $\vv$ and $\kk$, where each $k_i\geq 2$.  Then
\[ C(\vv,\kk,2) \leq C(\vmax,\kmin,2). \]
\end{thm}

\proof Using our results above, we obtain:
\[ \begin{array}{rcll}
C(\vv,\kk,2) & \leq & C(\vvmax,\kk,2)    & \textnormal{(by Lemma~\ref{lemma:deletion}, since $\vv \leq \vvmax$)} \\
             & \leq & C(\vvmax,\kkmin,2) & \textnormal{(by Lemma~\ref{lemma:expansion}, since $\kk \geq \kkmin$)} \\
             & =    & C(\vmax,\kmin,2)   & \textnormal{(by Corollary~\ref{cor:equiv}, with a single equivalence class)}
\end{array} \]
as required. \endproof

\begin{cor} \label{cor:minimaxpair}
Suppose we have vectors $\vv$ and $\kk$, with an index $i$ where $v_i=\vmax$ and $k_i=\kmin \geq 2$.  Then
\[ C(\vv,\kk,2) = C(\vmax,\kmin,2). \]
\end{cor}

\proof The inequality $C(\vv,\kk,2) \leq C(\vmax,\kmin,2)$ is given by Theorem~\ref{thm:minimax}.  We obtain reverse inequality by restricting to the part $i$ where $v_i=\vmax$ and $k_i=\kmin$, and applying Proposition~\ref{prop:restriction}.  \endproof

This last corollary is especially useful, as many of the covering numbers $C(v,k,2)$ are known exactly, and thus we are able to obtain the sizes of optimal generalized covering designs in many instances.  

\subsection{A construction algorithm} \label{subsect:construction}
All of our proofs in this section have been constructive.  As a direct consequence, we have an algorithm for actually constructing generalized covering designs from a single covering design.  Furthermore, if there is a part with $(v_i,k_i)=(\vmax,\kmin)$ and $\mathcal{C}$ is optimal, then (provided that $\kmin\geq 2$) Corollary~\ref{cor:minimaxpair} ensures that designs obtained using this construction are optimal.

\begin{construction} \label{constr:construction}
Suppose we are given vectors $\vv$ and $\kk$, where $k_i\geq 2$ for all $i$, and we have an optimal $(\vmax,\kmin,2)$-covering design $\mathcal{C}$.  Then we construct a ${\rm GC}(\vv,\kk,2)$ as follows:
\begin{itemize}
\item put a copy of $\mathcal{C}$ on each part;
\item in any part where $v_i < \vmax$, delete the extra points, replacing them with a ``placeholder'' symbol $\star$;
\item in any part where $k_i > \kmin$, add $k_i-\kmin$ placeholders to each block;
\item in each block, replace the placeholders greedily, ensuring that no symbol is repeated in a block;
\item remove any repeated blocks.
\end{itemize}
\end{construction}

We illustrate this construction with an example.

\begin{example} \label{example:fano}
Suppose $\vv=(5,6,7)$ and $\kk=(3,4,3)$.  So we have $\vmax=7$ and $\kmin=3$: an optimal $(7,3,2)$-covering design with 7 blocks is, of course, the Fano plane, whose blocks are
\[ \{124\}, \, \{235\}, \, \{346\}, \, \{457\}, \, \{156\}, \, \{267\}, \{137\}. \]

So our ``pre-design'' looks like:
\[ \begin{array}{ccc}
( \{124\},         & \{124\star\},     & \{124\} ) \\
( \{235\},         & \{235\star\},     & \{235\} ) \\
( \{34\star\},     & \{346\star\},     & \{346\} ) \\
( \{45\star\},     & \{45\star\star\}, & \{457\} ) \\
( \{15\star\},     & \{156\star\},     & \{156\} ) \\
( \{2\star\star\}, & \{26\star\star\}, & \{267\} ) \\
( \{13\star\},     & \{13\star\star\}, & \{137\} )
\end{array} \]
Having filled the placeholder positions, and using the least point available at each stage, we then obtain:
\[ \begin{array}{lll}
( \{124\}, & \{1234\}, & \{124\} ) \\
( \{235\}, & \{1235\}, & \{235\} ) \\
( \{134\}, & \{1346\}, & \{346\} ) \\
( \{145\}, & \{1245\}, & \{457\} ) \\
( \{125\}, & \{1256\}, & \{156\} ) \\
( \{123\}, & \{1236\}, & \{267\} ) \\
( \{123\}, & \{1234\}, & \{137\} )
\end{array} \]
As we have $(\vmax,\kmin)=(7,3)$ occurring in a part, we are guaranteed that this ${\rm GC}(\vv,\kk,2)$ is optimal.
\end{example}

In situations where there is no pair $(v_i,k_i)=(\vmax,\kmin)$, then it is possible to construct pathological examples which are far from optimal, even if beginning with an optimal covering design.  For instance, if $\vv=(100,7)$ and $\kk=(98,3)$, then this construction requires a $(100,3,2)$-covering design, which has size at least $C(100,3,2)\geq 1667$.  We will see later that 7 blocks can be used!

We conclude this subsection with a remark concerning the placeholder symbol $\star$.  In many constructions of covering arrays (see~\cite{Colbourn04,Colbourn06}, for instance), a ``don't care'' symbol is often appended to the alphabet.  This can then be replaced arbitrarily without affecting the requirement that all $t$-tuples of symbols be covered in every $t$-subset of columns, which can be particularly useful in applications.  So, in our construction, it may actually be beneficial to leave the placeholder symbols {\em in situ}, rather than filling those positions greedily.

\subsection{Amalgamation}
This is yet another operation to obtain a new generalized covering design from an existing one, this time by combining two parts into one.  

As usual, we let $\vv=(v_1,\ldots,v_m)$ and $\kk=(k_1,\ldots,k_m)$, and suppose $\mathcal{D}$ is a ${\rm GC(\vv,\kk,2)}$.  Suppose further that we have two indices where the entries of $\kk$ are at least~2, which without loss of generality we assume are $k_1$ and $k_2$.  Now let $\vv^+ = (v_1+v_2, v_3, \ldots, v_m)$ and $\kk^+ = (k_1+k_2,k_3,\ldots,k_m)$.  The operation of {\em amalgamation} allows us to construct a ${\rm GC}(\vv^+,\kk^+,2)$ from $\mathcal{D}$.  For each block $\BB=(B_1,\ldots,B_m) \in \mathcal{D}$, let $\BB^+ = (B_1 \dot\cup B_2, B_3, \ldots, B_m)$.  Finally, let $\mathcal{D}^+ = \{ \BB^+ \,\, : \,\, \BB\in\mathcal{D} \}$.

\begin{prop} \label{prop:amalg}
The design $\mathcal{D}^+$ defined above is a ${\rm GC}(\vv^+,\kk^+,2)$, and so $C(\vv^+,\kk^+,2) \leq C(\vv,\kk,2)$.
\end{prop}

\proof Once more, we think of the generalized covering design $\mathcal{D}$ as a clique covering of the graph $G_{\vv,\kk}$.  Now, observe that the graph $G_{\vv^+,\kk^+}$ must be isomorphic to $G_{\vv,\kk}$, except that the parts $X_1$ and $X_2$ of $V(G_{\vv,\kk})$, which have sizes $v_1$, $v_2$ respectively, can be thought of as one larger part of size $v_1+v_2$.  Each block $\BB\in\mathcal{D}$ corresponds to a clique in $G_{\vv,\kk}$ with $k_1$ vertices in $X_1$ and $k_2$ vertices in $X_2$, while the corresponding block $\BB^+\in\mathcal{D}^+$ gives exactly the same clique in $G_{\vv^+,\kk^+}$ (with $k_1+k_2$ vertices in $X_1 \dot\cup X_2$).  Since the two graphs are isomorphic, and $\mathcal{D}$ and $\mathcal{D}^+$ contain the same cliques, it follows that $\mathcal{D}^+$ covers all the edges in $G_{\vv^+,\kk^+}$.
\endproof 

The above result can be used to prove the following.

\begin{prop} \label{prop:amalg2}
Suppose we have vectors of positive integers $\vv=(v_1,\ldots,v_m)$ and $\kk=(k_1,\ldots,k_m)$ with $\vv\geq\kk$, and let $\rr=(r_1,\ldots,r_m)$ be a vector of integers satisfying $0\leq r_i \leq k_i-2$ for all $i$.  Then
\[ C(\vv,\kk,2) \leq C(\vv-\rr,\kk-\rr,2). \]
\end{prop}

\proof By repeatedly applying Proposition \ref{prop:amalg} above, we have
\[ C( (v_1,v_2,\ldots,v_m),(k_1,k_2,\ldots,k_m),2 ) \]
\[ \leq C((v_1-r_1,r_1,v_2-r_2,r_2,\ldots,v_m-r_m,r_m),(k_1-r_1,r_1,k_2-r_2,r_2,\ldots,k_m-r_m,r_m),2). \]
Then we apply Proposition \ref{prop:allvertexparts} to show that this is equal to
\[ C((v_1-r_1,v_2-r_2,\ldots,v_m-r_m),(k_1-r_1,k_2-r_2,\ldots,k_m-r_m),2) = C(\vv-\rr,\kk-\rr,2) \]
as required. \endproof

We remark that amalgamation can sometimes give much better results than the algorithm described in subsection~\ref{subsect:construction}.  Starting from the design in Example~\ref{example:fano}, we can use Proposition~\ref{prop:amalg} to obtain a ${\rm GC}((11,7),(7,3),2)$ with 7 blocks, by amalgamating the first and second parts.  This must be optimal, as it meets the bound in Corollary~\ref{cor:restriction}.  However, if we applied Construction~\ref{constr:construction} to obtain a ${\rm GC}((11,7),(7,3),2)$ from scratch, we would need to begin with an $(11,3,2)$-covering design.  The smallest such covering design has 19 blocks (see Gordon {\em et al.}\ \cite{GordonKuperbergPatashnik95}), meeting the Sch\"{o}nheim bound.  Taking that design, applying Construction~\ref{constr:construction} to it, and filling in the placeholders lexicographically, gives a ${\rm GC}((11,7),(7,3),2)$ with 19 blocks, which is considerably larger that that obtained using amalgamation.

Proposition \ref{prop:amalg2} can be combined with the argument of Theorem \ref{thm:minimax} to obtain the bound
\[ C(\vv,\kk,2) \leq C\left( \max_{1\leq j\leq m} \{v_j-(k_j-k_{\min})\},\, k_{\min},\, 2 \right) \]
(provided that $\kmin\geq 2$),which is often a considerable improvement.  For example, 
\[ C((100,7),(98,3),2)\leq C(\max\{100-(98-3),7-(3-3)\},3,2)=C(7,3,2), \]
which compares with the bound of $C(100,3,2)$ from Theorem~\ref{thm:minimax}.  (Note that $C(7,3,2)=7$, while $C(100,3,2)\geq 1667$, so this is definitely an improvement!)

\section{Another graphical interpretation}

There is an another interpretation of strength-2 generalized covering designs in terms of graphs, which  %%\marginpar{\sf [MODIFIED]}
while similar to that developed in Section 3, gives improved bounds on $C(\vv,\kk,2)$ for certain parameter sets.
Suppose we have vectors $\vv=(v_1,\ldots,v_m)$ and $\kk=(k_1,\ldots,k_m)$, with $\vv\geq\kk$, as usual.  As well as the interpretation in terms of clique coverings, we can also regard a ${\rm GC}(\vv,\kk,2)$ as an edge-covering of a complete multipartite graph, as explained below.

\begin{thm} \label{thm:GCgraph}
Suppose we have $\vv=(v_1,v_2,\ldots,v_m)$ and $\kk=(k_1,k_2,\ldots,k_m)$, where $\vv\geq\kk$.  Then a generalized covering design ${\rm GC}(\vv,\kk,2)$ is equivalent to an edge-covering of the complete multipartite graph $K_\vv$ by complete multipartite graphs $K_\kk$, where:
\begin{itemize}
\item[(i)] For each copy of $K_\kk$, the vertices contained in the part corresponding to $k_i$ are chosen from the vertices of $K_\vv$ corresponding to $v_i$;
\item[(ii)] the set of the complements of each copy of $K_\kk$ covers all edges of the complement of $K_\vv$ in parts where $k_i\geq 2$.
\end{itemize}
\end{thm}

\proof Each block of a ${\rm GC}(\vv,\kk,2)$ corresponds to a copy of the complete multipartite graph $K_\kk$ satisfying condition (i).  Now, we have two types of admissible vector $\tt$.  First, we have vectors $\tt$ consisting of two $1$s and the rest $0$s correspond to pairs of vertices in distinct parts: the fact that we have an edge-covering ensures that these pairs are contained within some block.  Second, we have vectors $\tt$ with exactly one $2$ in a position $i$ with $k_i\geq 2$, and $0$s elsewhere.  These correspond to pairs of vertices within a part, and condition (ii) ensures that these pairs are contained within some block. \endproof

We remark that in the case of covering arrays (i.e.\ where $\kk=(1,1,\ldots,1)$), the two interpretations are the same: in this case, the complete multipartite graph $K_\kk$ is actually a complete graph, so covering with copies of $K_\kk$ actually gives a clique covering.

%\noindent {\sf [BOUNDS OBTAINED FROM THIS]}
Using this second interpretation, we can obtain various bounds on $C(\vv,\kk,2)$, which do not always correspond to those obtained in Section~3.  The simplest of these is the na\"{\i}ve bound for covering a graph with isomorphic subgraphs, analogous to that in Proposition \ref{prop:naive}.  Simply by counting the number of edges in a complete multipartite graph, we arrive at the following result.

\begin{prop} \label{prop:naivemulti}
Where $\vv=(v_1,v_2,\ldots,v_m)$ and $\kk=(k_1,k_2,\ldots,k_m)$ and $\vv\geq\kk$, we have 
\[ C(\vv,\kk,2) \geq \left\lceil \frac{ \sum_{i\neq j} v_i v_j }{ \sum_{i\neq j} k_i k_j } \right\rceil. \]
\end{prop}

Depending on the precise nature of the entries of $\vv$ and $\kk$, it is quite possible that either of Propositions \ref{prop:naive} or \ref{prop:naivemulti} will give a better bound.

The notion of restriction, introduced in section 3.2, can also be applied to this interpretation.  Suppose $I$ is a subset of the index set $\{1,\ldots,m\}$.  Now, we can restrict $\vv$ and $\kk$ to $I$, and will still have a covering as described in Theorem \ref{thm:GCgraph}, provided that $|I|\geq 2$.  (If $|I|=1$, then both $K_{\vv^I}$ and $K_{\kk^I}$ have no edges, so the construction is vacuous.)  This gives rise to the following bound.

\begin{prop} \label{prop:restrictmulti}
For $\vv=(v_1,\ldots,v_m)$ and $\kk=(k_1,\ldots,k_m)$ where $\vv\geq\kk$, we have
\[ C(\vv,\kk,2) \geq \max_{\substack{ I \subseteq \{1,\ldots,m\} \\ |I|\geq 2}} \left\lceil \frac{|E(K_{\vv^I})|}{|E(K_{\kk^I})|} \right\rceil. \]
\end{prop}

As with the bound in Corollary \ref{cor:restriction}, this bound involves checking an exponential number of cases, so a more practical bound would be useful.  We can obtain one by restricting only to pairs of indices, i.e.\ where $I=\{i,j\}$ (for some $i,j\in\{1,\ldots,m\}$), as shown below.

\begin{prop} \label{prop:restrictpairs}
For $\vv=(v_1,\ldots,v_m)$ and $\kk=(k_1,\ldots,k_m)$ where $\vv\geq\kk$, we have
\[ C(\vv,\kk,2) \geq \max_{\stackrel{1 \leq i, j \leq m}{i\neq j}} \left\lceil \frac{v_i v_j}{k_i k_j} \right\rceil. \]
\end{prop}

%%\marginpar{\sf [MODIFIED]}
However, a refinement of this bound can be obtained as a special case of a bound on $C(\vv,\kk,t)$ for arbitrary values of $t$, as given in the next section.

%\begin{prop} \label{prop:restrictpairs2}
%For $\vv=(v_1,\ldots,v_m)$ and $\kk=(k_1,\ldots,k_m)$ where $\vv\geq\kk$, we have
%\[ C(\vv,\kk,2) \geq \max_{\stackrel{1 \leq i, j \leq m}{i\neq j}} \left\lceil \frac{v_i}{k_i} \left\lceil \frac{v_j}{k_j}\right\rceil\right\rceil. \]
%\end{prop}
%
%\proof Suppose $1 \leq i, j \leq m$ with $i \neq j$.  Corresponding to the vector $\tt$ with 1 in positions $i$ and $j$ and 0 elsewhere, we obtain that each pair $(x,y)$, where $x \in X_i$ and $y \in X_j$, appears in at least one block.  Let $x \in X_i$.  Then $x$ must occur in at least $\lceil v_j/k_j \rceil$ blocks, as otherwise, there would be an edge $xy$, where $y \in X_j$, which does not occur in any block.  Hence counting one for each block containing $x$, for each $x \in X_i$, gives at least $v_i \lceil v_j/k_j \rceil$.  However, since each block contains $k_i$ elements of $X_i$, we have counted each block $k_i$ times, and so the total number of blocks is at least
%\[ \frac{v_i}{k_i} \left\lceil \frac{v_j}{k_j} \right\rceil. \]
%The given bound follows since the number of blocks must be an integer.  \endproof

\section{A bound for arbitrary strength}

%%\marginpar{\sf [NEW\\ SECTION]}
When considering generalized covering designs of strength~$t\geq 2$, it is more difficult to obtain bounds on $C(\vv,\kk,t)$ as there are a variety of possible $(\kk,t)$-admissible vectors $\tt$.  However, one can obtain bounds by considering just one possible ``shape'' of such vectors.  In this section, we consider only the $(\kk,t)$-admissible vectors consisting of $t$ entries of~1 and all other entries~0, to obtain a bound which is an analogue of the Sch\"onheim bound (Theorem~\ref{thm:Schonheim}).

\begin{prop} \label{prop:revised_bound}
Suppose that $\vv=(v_1,\ldots,v_m)$ and $\kk=(k_1,\ldots,k_m)$ where $\vv\geq\kk$, and that $t\leq m$.  Let $\{i_1, i_2, \ldots, i_t\} \subseteq \{1, 2, \ldots, m\}$, and let $\mathcal{C}$ be a collection of blocks in ${\XX \choose \kk}$ which contain each $t$-tuple of the form $(x_{i_1}, x_{i_2}, \ldots, x_{i_t})$, where $x_{i_j} \in X_{i_j}$.  Then 
\[
|\mathcal{C}| \geq \left\lceil \frac{v_{i_1}}{k_{i_1}} \left\lceil \frac{v_{i_2}}{k_{i_2}} \cdots \left\lceil \frac{v_{i_t}}{k_{i_t}} \right\rceil \cdots \right\rceil \right\rceil.
\]
\end{prop}

\begin{proof}
The proof is by induction on $t$.  If $t=1$, then the result is obvious.  

%For the induction hypothesis, we 
Suppose that the result holds for %(any $v, k, m \geq t$ and) 
any set $\{j_1, j_2, \ldots, j_{t-1}\} \subseteq \{1,2, \ldots, m\}$.  
Now choose $\{i_1, i_2, \ldots, i_t\} \subseteq \{1, 2, \ldots, m\}$, and let $\mathcal{C}$ be as defined in the statement of the proposition.  
%Note that the number of elements of $X_{i_1}$ (with repetition) in $\mathcal{C}$ is $k_{i_1} |\mathcal{C}|$.
Note that the total number of occurences (with repetition) of elements of $X_{i_1}$ in $\mathcal{C}$ is $k_{i_1} |\mathcal{C}|$.  
Moreover, for each element $x \in X_{i_1}$, $\mathcal{C}$ must contain the $t$-tuple $(x, x_{i_2}, \ldots, x_{i_t})$ for any choice of elements $x_{i_j} \in X_{i_j}$, where $j=2, \ldots, t$.  By the induction hypothesis, it follows that $\mathcal{C}$ must contain at least
\[
\left\lceil \frac{v_{i_2}}{k_{i_2}} \left\lceil \frac{v_{i_3}}{k_{i_3}} \cdots \left\lceil \frac{v_{i_t}}{k_{i_t}} \right\rceil \cdots \right\rceil \right\rceil
\]
blocks for each element of $X_{i_1}$.  In total, this means that the 
number of occurences (with repetition) of elements of $X_{i_1}$ in $\mathcal{C}$ is at least
\[
v_{i_1} \left\lceil \frac{v_{i_2}}{k_{i_2}} \left\lceil \frac{v_{i_3}}{k_{i_3}} \cdots \left\lceil \frac{v_{i_t}}{k_{i_t}} \right\rceil \cdots \right\rceil \right\rceil.
\]
Hence 
\[
k_{i_1} |\mathcal{C}| \geq v_{i_1} \left\lceil \frac{v_{i_2}}{k_{i_2}} \left\lceil \frac{v_{i_3}}{k_{i_3}} \cdots \left\lceil \frac{v_{i_t}}{k_{i_t}} \right\rceil \cdots \right\rceil \right\rceil,
\]
and so we conclude that 
\[
|\mathcal{C}| \geq \left\lceil \frac{v_{i_1}}{k_{i_1}} \left\lceil \frac{v_{i_2}}{k_{i_2}} \cdots \left\lceil \frac{v_{i_t}}{k_{i_t}} \right\rceil \cdots \right\rceil \right\rceil.
\]
\end{proof}

By considering all possible choices of $\{i_1, \ldots i_t\}$, we obtain the following bound as a direct consequence.

\begin{cor} \label{cor:allposs}
If $\vv=(v_1,\ldots,v_m)$ and $\kk=(k_1,\ldots,k_m)$ where $\vv\geq\kk$, and if $t \leq m$, then 
\[
C(\vv, \kk, t) \geq \max_{\{i_1, \ldots i_t\} \subseteq\{1, \ldots, m\}} \left\lceil \frac{v_{i_1}}{k_{i_1}} \left\lceil \frac{v_{i_2}}{k_{i_2}} \cdots \left\lceil \frac{v_{i_t}}{k_{i_t}} \right\rceil \cdots \right\rceil \right\rceil. 
\]
\end{cor}

\begin{proof}
Given $\{i_1, \ldots i_t\} \subseteq\{1, \ldots, m\}$, consider the vector $\tt$ which has 1 in positions $i_1, \ldots, i_t$, and 0 elsewhere; this vector will always be $(\kk,t)$-admissible regardless of $\kk$.  Corresponding to the vector $\tt$, we obtain that the set $\mathcal{B}$ of blocks of a $\mathrm{GC}(\vv,\kk,t)$ must contain each $t$-tuple of the form $(x_{i_1}, x_{i_2}, \ldots, x_{i_t})$, where $x_{i_j} \in X_{i_j}$.  By Proposition~\ref{prop:revised_bound}, the number of blocks must therefore be at least 
\[ \left\lceil \frac{v_{i_1}}{k_{i_1}} \left\lceil \frac{v_{i_2}}{k_{i_2}} \cdots \left\lceil \frac{v_{i_t}}{k_{i_t}} \right\rceil \cdots \right\rceil \right\rceil. \]
\end{proof}

In the special case where $t=2$, this reduces to the following. 

\begin{cor} \label{cor:restrictpairs2}
%%\marginpar{\sf [WAS PROP. 4.5]}
For $\vv=(v_1,\ldots,v_m)$ and $\kk=(k_1,\ldots,k_m)$ where $\vv\geq\kk$, we have
\[ C(\vv,\kk,2) \geq \max_{\stackrel{1 \leq i, j \leq m}{i\neq j}} \left\lceil \frac{v_i}{k_i} \left\lceil \frac{v_j}{k_j}\right\rceil\right\rceil. \]
\end{cor}

We remark that this second corollary gives the desired refinement of the bound in Proposition~\ref{prop:restrictpairs}.

\section{Product constructions}

\subsection{The block-recursive construction for arbitrary strength} \label{subsect:blockrec}

Our first construction is based on the block-recursive construction for covering arrays that appears in Poljak {\em et al.}~\cite{Poljak83} and also Stevens and Mendelsohn~\cite{Stevens97}, and which was later extended to mixed covering arrays by Colbourn {\em et al.}~\cite{Colbourn06}.  This construction uses two strength-$t$ generalized covering designs to construct another strength-$t$ generalized covering design.  We recall the notion of concatenation of vectors from Section \ref{section:defn}.

%To begin, we need to define a simple operation on $m$-tuples.
%Let $\vv = (v_1,v_2,\ldots,v_m)$ and $\ww = (w_1,w_2,\ldots,w_n)$ be an
%$m$-tuple and an $n$-tuple respectively. Define the {\em concatenation} of $\vv$
%and $\ww$ to be the $(m+n)$-tuple
%\[
%\cat(\vv,\ww) = (v_1,v_2,\ldots,v_m,w_1,w_2,\ldots,w_n).
%\] 

\begin{thm} \label{thm:blockrec}
  Suppose $\vv=(v_1,v_2,\ldots,v_m)$ and $\kk=(k_1,k_2,\ldots,k_m)$
  are $n$-tuples with $\vv\geq\kk$ and $\ww=(w_1,w_2,\ldots,w_n)$
  and $\bl=(\ell_1,\ell_2,\ldots,\ell_n)$ are $n$-tuples with $\ww\geq\bl$.  
  Assume that $\mathcal{D}_1$ is a generalized covering design ${\rm GC}(\vv,\kk,t)$ with $b$ blocks, and that $\mathcal{D}_2$ is a ${\rm GC}(\ww,\bl,t)$ with $c$ blocks. Then there exists a ${\rm GC}(\cat(\vv,\ww), \cat(\kk,\bl), t)$ with $bc$ blocks.
\end{thm}

\proof Assume that $\mathcal{B}=\{\BB_1,\BB_2,\ldots, \BB_b\}$ are the blocks
of $\mathcal{D}_1$ and $\mathcal{C} = \{\CC_1,\CC_2,\ldots,\CC_c\}$ are the
blocks of $\mathcal{D}_2$.

Since each $\BB_i$ is an $m$-tuple of sets of sizes $(k_1,k_2,\ldots,k_m)$, and each $\CC_j$
is an $n$-tuple of sets of sizes $(\ell_1,\ell_2,\ldots,\ell_n)$, we have that $\cat(\BB_i,\CC_j)$ is
the $(m+n)$-tuple of sets (of sizes labelled by $\cat(\kk,\bl)$) formed by
the concatenation of $\BB_i$ and $\CC_j$. 

We claim that the set of blocks
\[ \left\{ \cat(\BB_i,\CC_j) \; : \; i = 1,2,\ldots,b, \; j=1,2,\ldots,c \right\} \]
form a generalized covering design ${\rm GC}(\cat(\vv,\ww), \cat(\kk,\bl),t)$ with $bc$ blocks.

To see that this claim is true, let $\TT =(T_1,T_2,\ldots,T_{m+n})$ be
any $(m+n)$-tuple of sets that is
$(\cat(\vv,\ww),\cat(\kk,\bl),t)$-admissible. We will prove that there
exists some block $\cat(\BB_i,\CC_j)$ that contains $\TT$.

Consider the $m$-tuple $\TT' = (T_1,T_2,\ldots,T_m)$ formed by taking the
first $m$ sets in $\TT$.  Clearly, $\sum_{i=1}^n |T_i| \leq t$, and
since $\mathcal{D}_1$ is a generalized covering design, by
Proposition~\ref{prop:reducet}, the $n$-tuple $\TT'$ is contained in some block
$\BB_i$ of $\mathcal{B}$.

Similarly, if we define $\TT''= (T_{m+1}, T_{m+2}, \ldots, T_{m+n})$, there
is a block $\CC_{j}$ of $\mathcal{D}_2$ that contains $\TT''$.  Thus
$\TT = \cat(\TT',\TT'')$ is contained in the block $\cat(\BB_i,\CC_j)$.  \endproof

This construction can be used to get an upper bound on the size of generalized covering designs.

\begin{cor}\label{cor:bcbound}
  Suppose we have $\vv=(v_1,v_2,\ldots,v_m)$ and $\kk=(k_1,k_2,\ldots,k_m)$
  with $\vv\geq\kk$, and $\ww=(w_1,w_2,\ldots,w_m)$ and
  $\bl=(\ell_1,\ell_2,\ldots,\ell_m)$ with $\ww\geq\bl$.  Then for
  all appropriate values of $t$,
\[
C(\cat(\vv,\ww), \cat(\kk,\bl), t) \leq C(\vv,\kk,t) \; C(\ww,\bl,t). 
\]
\end{cor}

Unfortunately, this construction can lead to very poor upper
bounds on the number of blocks in a generalized covering design.  For example,
consider a ${\rm GC}(\vv,\kk,2)$ with the condition that no entries of $\kk$
are equal to one.  Then the block-recursive construction produces the bound
\[
C(\cat(\vv,\vv), \cat(\kk,\kk), 2) \leq \left( C(\vv,\kk,2) \right)^2, 
\]
whereas it follows from Corollary~\ref{cor:equiv} the actual value of $C(\cat(\vv,\vv),\cat(\kk,\kk), 2)$ is equal to $C(\vv, \kk, 2)$ (again, provided that all of the entries of $\kk$ are at least two).

However, there is a modification of this construction that in some
circumstances produces a better bound than the one given in
Corollary~\ref{cor:bcbound}. Again, suppose that
$\vv=(v_1,v_2,\ldots,v_m)$ and $\kk=(k_1,k_2,\ldots,k_m)$ with $\vv\geq\kk$, that
$\ww=(w_1,w_2,\ldots,w_n)$ and $\bl=(\ell_1,\ell_2,\ldots,\ell_n)$ with $\ww\geq\bl$, 
and that we have two generalized covering designs, a ${\rm GC}(\vv,\kk,t)$ and a ${\rm GC}(\ww,\bl,t)$. 
Assume that $\mathcal{B}=\{\BB_1,\BB_2,\ldots, \BB_b\}$ are the blocks of the first design 
and $\mathcal{C} = \{\CC_1,\CC_2,\ldots,\CC_c\}$ are the blocks of the second design, 
and further assume that $b\geq c$.

The key to this improvement is to consider the set of $b$ blocks
\[ \mathcal{S} = \{ \cat(\BB_1,\CC_1) , \cat(\BB_2,\CC_2), \; \dots \; \cat(\BB_c, \CC_c) \} \cup \{ \cat(\BB_{c+1}, \CC_1),\; \dots \; \cat(\BB_b,\CC_1) \}. \]
(In the last $b-c$ blocks of~$\mathcal{S}$, the $m$-tuple $\CC_1$ could be replaced with
any block of $\mathcal{C}$.)  We note that $\mathcal{S}$ contains any $(\cat(\vv,\ww),\cat(\kk,\bl),t)$-admissible vector $\TT =(T_1,T_2,\ldots,T_{m+n})$ with the condition that either 
\[
\sum_{i=1}^m |T_i| = 0 \quad  \textrm{or} \quad \sum_{i=m+1}^{m+n} |T_i| = 0.
\]

To extend the set~$\mathcal{S}$ to a generalized covering design, we need to
add blocks that will cover all the remaining
$(\cat(\vv,\ww),\cat(\kk,\bl),t)$-admissible vectors. These vectors
will all be of the form $\TT =(T_1,T_2,\ldots,T_{m+n})$ where
\[
\sum_{i=1}^m |T_i| \leq t-1 \quad  \textrm{and} \quad \sum_{i=m+1}^{m+n} |T_i| \leq t-1.
\]
It is possible to construct such blocks by applying the
block-recursive construction to a \mbox{${\rm GC}(\vv,\kk,t-1)$} and a
${\rm GC}(\ww,\bl,t-1)$.  Using this modification to the block-recursive
construction, we obtain the following bound.

\begin{thm} \label{thm:blockrecsimp}
  Suppose we have $\vv=(v_1,v_2,\ldots,v_m)$ and $\kk=(k_1,k_2,\ldots,k_m)$ with $\vv\geq\kk$, 
and $\ww=(w_1,w_2,\ldots,w_n)$ and $\bl=(\ell_1,\ell_2,\ldots,\ell_n)$ with $\ww\geq\bl$. 
Then for all suitable $t>0$,
\[
C(\cat(\vv,\ww), \cat(\kk,\bl), t) \leq \max\{ C(\vv,\kk,t),C(\ww,\bl,t)\} 
   + C(\vv,\kk,t-1) \; C(\ww,\bl,t-1). 
\]
\end{thm}

Two special cases are worth mentioning.  First, in the case $t=1$, we obtain
\[ C(\cat(\vv,\ww), \cat(\kk,\bl), 1) \leq \max\{ C(\vv,\kk,1),C(\ww,\bl,1)\}, \]
(thanks to the convention that $C(\vv,\kk,0)=0$); in fact, Proposition~\ref{prop:tis1} ensures that this bound holds with equality.
Second, in the case where $t=2$, Theorem~\ref{thm:blockrecsimp} can be simplified by applying Proposition~\ref{prop:tis1}.

\begin{cor} \label{cor:blockrecsimp2}
  Suppose we have $\vv=(v_1,v_2,\ldots,v_m)$ and $\kk=(k_1,k_2,\ldots,k_m)$ with $\vv\geq\kk$, 
and $\ww=(w_1,w_2,\ldots,w_n)$ and $\bl=(\ell_1,\ell_2,\ldots,\ell_n)$ with $\ww\geq\bl$.  
Then the size of a ${\rm GC}(\cat(\vv,\ww), \cat(\kk,\bl), 2)$ is no more than
\[
\max \left\{ C(\vv,\kk,2), \; C(\ww,\bl,2) \right\} + \left( \max_{i=1,\ldots,n} \ceilfrac{v_i}{k_i} \right) \left( \max_{i=1,\ldots,m} \ceilfrac{w_i}{\ell_i} \right).
\]
\end{cor}

Note that, because the modified block-recursive construction of
Theorem~\ref{thm:blockrecsimp} has two stages, it is possible that the
design it produces may have repeated blocks (even if the inputted
designs did not).  We now illustrate this construction of
Theorem~\ref{thm:blockrecsimp} with an example (and one where that
possibility arises).

\begin{example}
Let $\BB_1,\ldots,\BB_{10}$ below be the blocks of a ${\rm GC}((5,7),(2,3), 2)$.
\[
\begin{array}{lll}
\BB_1= &  ( \{12\},  & \{123\} ) \\
\BB_2= &  ( \{34\}, & \{147\} ) \\
\BB_3= &  ( \{15\}, & \{156\} ) \\
\BB_4= &  ( \{45\}, & \{246\} ) \\
\BB_5= &  ( \{23\}, & \{257 \} )\\
\BB_6= &  ( \{24\}, & \{345\} ) \\
\BB_7= &  ( \{35\}, & \{367\} ) \\
\BB_8= &  ( \{14\}, & \{124\} ) \\
\BB_9= &  ( \{13\}, & \{127\} ) \\
\BB_{10}= &  ( \{25\}, &  \{126\} )\\
\end{array}
\]

Further, let $\CC_1,\ldots,\CC_{6}$ below be the blocks of a ${\rm GC}((3,4) , (2,2), 2)$.
\[ 
\begin{array}{lll}
\CC_1= &  ( \{12\}, & \{12\} )\\
\CC_2= &  ( \{13\}, & \{13\} )\\
\CC_3= &  ( \{12\}, & \{14\} )\\
\CC_4= &  ( \{23\}, & \{23\} )\\
\CC_5= &  ( \{23\}, & \{24\} )\\
\CC_6= &  ( \{13\}, & \{34\} )\\
\end{array}
\]

These two designs can be used to build a ${\rm GC}((5,7,3,4), (2,3,2,2), 2)$
with the $16$ blocks $\DD_1,\ldots,\DD_{16}$ given below.
\[
\begin{array}{lllll}
\DD_1= &  ( \{12 \}, & \{123 \}, & \{12 \}, & \{12 \}) \\
\DD_2= &  ( \{34 \}, & \{147 \}, & \{13 \}, & \{13 \}) \\
\DD_3= &  ( \{15 \}, & \{156 \}, & \{12 \}, & \{14 \}) \\
\DD_4= &  ( \{45 \}, & \{246 \}, & \{23 \}, & \{23 \}) \\
\DD_5= &  ( \{23 \}, & \{257 \}, & \{23 \}, & \{24 \}) \\
\DD_6= &  ( \{24 \}, & \{345 \}, & \{13 \}, & \{34 \}) \\
\DD_7= &  ( \{35 \}, & \{367 \}, & \{12 \}, & \{12 \}) \\
\DD_8= &  ( \{14 \}, & \{124 \}, & \{12 \}, & \{12 \}) \\
\DD_9= &  ( \{13 \}, & \{127 \}, & \{12 \}, & \{12 \}) \\
\DD_{10}= &  (\{25 \}, & \{126 \}, & \{12 \}, &  \{12 \}) \\
\DD_{11}= &  ( \{12 \}, &  \{123 \}, &\{12 \}, & \{12 \}) \\
\DD_{12}= &  ( \{34 \}, &  \{456 \}, &\{12 \}, & \{12 \}) \\
\DD_{13}= &  ( \{15 \}, &  \{127 \}, &\{12 \}, & \{12 \}) \\
\DD_{14}= &  ( \{12 \}, &  \{123 \}, &\{13 \}, & \{34 \}) \\
\DD_{15}= &  ( \{34 \}, &  \{456 \}, &\{13 \}, & \{34 \}) \\
\DD_{16}= &  ( \{15 \}, &  \{127 \}, &\{13 \}, & \{34 \}) \\
\end{array}
\]

The blocks $\DD_1$ through to $\DD_{10}$ are formed as the concatenation of 
blocks of the original designs. The final $6$ blocks are formed by
applying the block-recursive construction on a ${\rm GC}( (5,7),(2,3),1)$ and a
${\rm GC}( (3,4),(2,2),1)$.

Note that $\DD_1 = \DD_{11}$, so we have a repeated block, and thus $\DD_{11}$ may be deleted leaving a $\rm{GC}((3,4,5,7), (2,2,2,3), 2)$ with 15 blocks.  We note that Theorem~\ref{thm:minimax} gives an upper bound of $C(7,2,2)=21$ blocks for a generalized covering design with these parameters, while our various lower bounds show that the minimum number of blocks is at least 10 (this is obtained by using Corollary~\ref{cor:restriction2}).
\end{example}

\subsection{A MacNeish-type construction for strength 2} \label{subsect:MacNeish}
In a 1922 paper~\cite{MacNeish22}, MacNeish gave a recursive construction for mutually orthogonal Latin squares of order $mn$ from those of orders $m$ and $n$.  This was subsequently generalized to orthogonal arrays in 1952 by Bush~\cite{Bush52}, whose construction can also be applied to covering arrays.  Our second recursive construction for generalized covering designs is based on MacNeish's approach.

This construction takes two strength-$2$ generalized covering designs and produces a
third strength-$2$ generalized covering design.  As was the case with the block-recursive construction, we need to define an operation on
vectors to define this construction.

Let $\vv = (v_1,v_2,\ldots,v_m)$ and $\ww = (w_1, w_2, \ldots, w_m)$ be two
$m$-tuples of positive integers. Define the {\em Hadamard product} of $\vv$
and $\ww$ to be
\[
\vv \circ \ww = (v_1w_1, v_2w_2,\ldots,v_mw_m).
\]

We also need to define an operation on sets. Let
$R=\{r_1,r_2,\ldots,r_k\}$ be a subset of $\{1,2,\ldots,v\}$ and $S = \{s_1,
s_2, \ldots, s_\ell\}$ a subset of $\{1,2,\ldots,w\}$.  Then define
\[
R \circ_v S = \{r_i + (s_j - 1) v \; | \; i = 1,2,\ldots,k,\; j = 1,2,\ldots,\ell \},
\]
noting that $T \circ_v S$ is a subset of $\{1,\ldots, vw\}$ of size $k\ell$.  Now, if
$\RR = (R_1,R_2,\ldots,R_m)$ and $\SS = (S_1,S_2,\ldots,S_m)$ are two
$m$-tuples of sets and $\vv = (v_1,v_2,\ldots,v_m)$, we define
\[
\RR \circ_{\vv} \SS = (R_1 \circ_{v_1} S_1, R_2 \circ_{v_2} S_2, \ldots, R_m \circ_{v_m} S_m).
\]

\begin{thm} \label{thm:MacNeish}
Suppose that $\vv=(v_1,v_2,\ldots,v_m)$ and $\kk=(k_1,k_2,\ldots,k_m)$ with $\vv\geq\kk$, and that $\ww=(w_1,w_2,\ldots,w_m)$ and $\bl=(\ell_1,\ell_2,\ldots,\ell_m)$ with $\ww\geq\bl$.  Now suppose we have generalized covering designs ${\rm GC}(\vv,\kk,2)$ and ${\rm GC}(\ww,\bl,2)$ 
with $b$ and $c$ blocks respectively.  Then there exists a \mbox{${\rm GC}( (\vv \circ \ww), (\kk \circ \bl), 2)$} with $bc$ blocks.
\end{thm}

\proof
Let $\mathcal{B}=\{\BB_1,\BB_2,\ldots, \BB_b\}$ be the blocks of the first design
and $\mathcal{C} = \{\CC_1,\CC_2,\ldots,\CC_c\}$ be the blocks of the second
design.

Then we claim that the collection 
\[ \mathcal{D} = \{ \BB_i \circ_{\vv} \CC_j \; : \; i=1,\ldots,b, \; j = 1,\ldots,c \} \]
of $m$-tuples of sets forms the blocks of a ${\rm GC}( (\vv \circ \ww), (\kk \circ \bl), 2)$.  Clearly, $|\mathcal{D}|=bc$.

Let $\TT = (T_1,T_2,\ldots,T_m)$ be any $( (\vv \circ \ww), (\kk \circ \bl),2)$-admissible vector. We will show that this is contained in
some block of the form $\BB_i \circ_{\vv} \CC_j$.

First, consider the case where $\TT$ has one entry a set of size~2, and all other entries are the empty set.
In particular, we have that $T_i= \{t_1,t_2\}$ for some $i \in \{1,\ldots,m\}$ and that $T_j =
\emptyset$ for all $j\neq i$.  Then the set $T_i' = \{ t_1 \!\!\pmod{v_i}, t_2 \!\!\pmod{v_i} \}$ is a subset of $\{1,2,\ldots,v_i\}$.  (We use the convention that if $v_i\mid t_j$, then $t_j \!\!\pmod{v_i}=v_i$, rather than 0.)
%this should be $X_i$ but I don't have this defined!!
Construct an $m$-tuple $\TT'$ of sets by taking the $i^{th}$ entry of $\TT'$ to be $T_i'$ and all other entries the empty set. Then $\TT'$
is a $(\vv, \kk, 2)$-admissible vector. Since $\mathcal{B}$ is the set of blocks of a generalized covering design, there is a block $\BB_i$ that contains $\TT'$.

Similarly, the set $T_i''=\left\{ \ceilfrac{t_1}{v_i}, \ceilfrac{t_2}{v_i}\right\}$ is a
subset of $\{1,2,\ldots,w_i\}$ (since both $t_1\leq v_iw_i$ and $t_2\leq v_iw_i$). Then the $m$-tuple $\TT''$ that has $T_i''$ in the
entry $i$ and the empty set in all other entries is a $(\ww,\bl,2)$-admissible vector. Since $\mathcal{C}$ is the set of blocks of
a generalized covering design, there is a block $\CC_j$ that contains $\TT''$.

Now we notice that the block $\BB_i \circ_{\vv} \CC_j$ contains $\TT$, since
\[
t_1= t_1 \!\!\pmod{v_i} + \left (\ceilfrac{t_1}{v_i} - 1\right) v_i  \,\, \textrm{and} \\,
t_2 = t_2 \!\!\pmod{v_i} + \left(\ceilfrac{t_2}{v_i} - 1 \right) v_i.
\]

A similar argument also works for the case where the vector $\TT$ has two entries which are singletons and all other entries are the empty set.  To see this, assume that $\TT = (T_1, T_2, \ldots, T_m)$ with $T_i= \{t_1\}$ and $T_j= \{t_2\}$ (where $i \neq j$), and that $T_k = \emptyset$ for all $k \in
\{1,2,\ldots,n\} \setminus \{i,j\}$.
 
Consider the $n$-tuples
$\TT' =  (T_1', T_2', \ldots, T_n')$ where
\[
T_i'= \{t_1 \!\!\pmod {v_i} \}, \; T_j' = \{t_2 \!\!\pmod{v_j} \}, \; T_k' = \emptyset \; \textrm{for $k \neq i,j$},
\]
and
$\TT'' =  (T_1'', T_2'', \ldots, T_n'')$ where
\[
T_i''= \left\{ \ceilfrac{t_1}{v_i} \right\}, \; T_j'' = \left\{\ceilfrac{t_2}{v_j} \right\}, \;
T_k'' = \emptyset \quad \textrm{for $k \neq i,j$.}
\]
As before, there is a block $\BB_i \in \mathcal{B}$ that contains $\TT'$
and a block $\CC_j$ of $\mathcal{C}$ that contains $\TT''$.  Consequently, the block
$\BB_i \circ_\vv \CC_j$ contains the vector $\TT$. 
\endproof

Once again, we illustrate this construction with an example.

\begin{example}
Below is an example of a ${\rm GC}((3,4) , (2,3), 2)$ with 3 blocks:
\[
\begin{array}{lll}
\BB_1=&  ( \{12\}, & \{123\} ) \\
\BB_2=&  ( \{13\}, & \{124\} ) \\
\BB_3=&  ( \{23\}, & \{134\} ) \\
\end{array}
\]

From this we can construct the ${\rm GC}((9,16), (4,9), 2)$ with $9$ blocks given below:
% \[
% \begin{array}{c|c}
% \{1,2\} \circ \{1,2\} & \{1,2,3\}\circ \{1,2,3\} \\
% \{1,2\} \circ \{1,3\} & \{1,2,3\}\circ \{1,2,4\} \\
% \{1,2\} \circ \{2,3\} & \{1,2,3\}\circ \{1,3,4\} \\
% \{1,3\} \circ \{1,2\} & \{1,2,4\}\circ \{1,2,3\} \\
% \{1,3\} \circ \{1,3\} & \{1,2,4\}\circ \{1,2,4\} \\
% \{1,3\} \circ \{2,3\} & \{1,2,4\}\circ \{1,3,4\} \\
% \{2,3\} \circ \{1,2\} & \{1,3,4\}\circ \{1,2,3\} \\
% \{2,3\} \circ \{1,3\} & \{1,3,4\}\circ \{1,2,4\} \\
% \{2,3\} \circ \{2,3\} & \{1,3,4\}\circ \{1,3,4\} \\
% \end{array}
% \]

\[
\begin{array}{lll}
\DD_1 = &(\{1,2,4,5\}, & \{1,2,3,5,6,7,9,10,11\} )\\
\DD_2 = &(\{1,2,7,8\}, & \{1,2,3,5,6,7,13,14,15\} )\\
\DD_3 = &(\{4,5,7,8\}, & \{1,2,3,9,10,11,13,14,15\} )\\
\DD_4 = &(\{1,3,4,6\}, & \{1,2,4,5,6,8,9,10,12\})\\
\DD_5 = &(\{1,3,7,9\}, & \{1,2,4,5,6,8,13,14,16\})\\
\DD_6 = &(\{4,6,7,9\}, & \{1,2,4,9,10,12,13,14,16\})\\
\DD_7 = &(\{2,3,5,6\}, & \{1,3,4,5,7,8,9,11,12\})\\
\DD_8 = &(\{2,3,8,9\}, & \{1,3,4,5,7,8,13,15,16\})\\
\DD_9 = &(\{5,6,8,9\}, & \{1,3,4,9,11,12,13,15,16\})
\end{array}
\]

\end{example}

\section{Conclusion}
%{\sf [DIRECTIONS FOR FUTURE WORK]} \\ 
In this paper, we have defined a very broad class of objects.  Consequently, there are a considerable variety of directions for future work on this topic.  Our work in this paper has primarily been focused on strength-$2$ generalized covering designs, principally because of the interpretation in terms of clique coverings.  Two natural directions would be (i) to obtain further results in this case, and (ii) to extend our results to values of $t\geq 2$.

There is an improvement on the block-recursive construction for covering arrays that uses two copies of the same covering array with {\em disjoint blocks} removed from the second array~\cite{Stevens97}. With this improvement, covering arrays can be constructed that meet the asymptotic bound due to Gargano {\em et al.}~\cite{Gargano92}. It would be interesting to see if this improvement could also be applied to the block-recursive construction for generalized covering designs in subsection~\ref{subsect:blockrec}.

Wilson gave a recursive construction for transversal designs (described in~\cite{Beth99}) that has been applied to both covering arrays and covering arrays with mixed alphabets. (MacNeish's construction is a special case of this.)  It would be worthwhile to investigate how Wilson's construction can be adapted to generalized covering designs.

If we were to consider generalized covering designs of higher strength (i.e.\ with $t\geq 2$), the natural extension of our approach would be to define a notion of clique coverings in $t$-uniform hypergraphs.  Thus if a clique in a $t$-uniform hypergraph is defined to be a set $K$ of vertices
such that any subset of $K$ of size at most $t$ is contained in a hyperedge, then a strength-$t$ generalized covering design would be equivalent to a
clique covering of an appropriately-constructed $t$-uniform hypergraph.

Another direction for further work on generalized covering designs is to conduct an analysis of families of these designs with specific parameters. For example, one family that could be examined is that of all generalized covering designs with $\vv = (v_1,v_2)$.  Such a design would be, in some sense, close to a covering design. Alternatively, one could consider the all generalized covering designs for which $\kk = (x,1,1,\dots,1)$ with $x>1$; such designs would not be far removed from covering arrays.  

For small values of $k$ one could consider all generalized covering designs with $\kk = (k_1,k_2,\ldots,k_m)$ with $\sum_i k_i = k$.  In~\cite{Cameron09}, Cameron gives a complete description of generalized $t$-designs for $k\leq 4$ and $t\leq 3$: such a characterization could potentially be obtained for generalized covering designs.  With small values of $v =\sum v_i$, computer searches could be implemented to compare our various constructions.  With sufficiently many parameters fixed, it may also be possible to determine the asymptotic growth of the optimal size of such a generalized covering design.

In \cite{Cameron09}, Cameron also suggests the dual notion of generalized packing designs; the first two authors have also been investigating this problem~\cite{packings}.  While the definition is very similar to that of generalized covering designs, as is often the case with packing problems, the theory turns out to be quite different.  Various previously-known classes of designs arise as generalized packing designs; examples include Howell designs, Room squares, Hanani triple systems, and mutually orthogonal Latin rectangles. %%\marginpar{\sf [ADDITION]}

Finally, we recall that much of the motivation for work on covering designs and covering arrays was their widespread use in applications (particularly in communications and software testing).  It would be very interesting to discover applications for other classes of generalized covering designs.

\section*{Acknowledgements}
%UofR DMRG, PIMS, NSERC
The work in this paper was a joint project of the Discrete Mathematics Research Group at the University of Regina, attended by all the authors.  They would like to thank the other members of the group (B.~Ahmadi, F.~Alinaghipour, S.~M.~Fallat, J.~C.~Fisher, A.~M.~Purdy and S.~Zilles) for their participation.  The authors also thank the two anonymous referees for their helpful suggestions.  R.~F.~Bailey is a PIMS Postdoctoral Fellow.  A.~C.~Burgess and M.~S.~Cavers are NSERC Postdoctoral Fellows.  K.~Meagher acknowledges support from an NSERC Discovery Grant.

\end{document}